\documentclass[11pt]{article}
\usepackage{epsfig}
\usepackage{psfrag}
\usepackage{graphicx}
\usepackage{float,enumitem}
\usepackage[latin1]{inputenc}
\usepackage[colorlinks=true]{hyperref}
\usepackage{latexsym}
\usepackage{amssymb}
\usepackage{fullpage}
\def\Prob{{\hbox{\rm Prob}}}

\def\conv{{\rm Conv \,}}

\def\nr{\par \noindent}

\def\Def{\stackrel{\mathrm{def}}{=}}

\def\dom{{\rm dom \,}}
\def\beq{\begin{equation}}
\def\eeq{\end{equation}}

\newtheorem{theorem}{Theorem}
\newtheorem{lemma}{Lemma}
\newtheorem{corollary}{Corollary}

\newtheorem{assumption}{Assumption}
\newtheorem{definition}{Definition}

\newtheorem{example}{Example}
\newtheorem{remark}{Remark}
\newcommand{\proof}{\bf Proof: \rm \nr}
\newcommand{\half}{\mbox{\small ${1 \over 2}$}}

\def\ba{\begin{array}}
\def\ea{\end{array}}
\def\beann{\begin{eqnarray*}}
\def\eeann{\end{eqnarray*}}
\def\bea{\begin{eqnarray}}
\def\eea{\end{eqnarray}}

\def\BT{\begin{theorem}}
\def\ET{\end{theorem}}
\def\BL{\begin{lemma}}
\def\EL{\end{lemma}}
\def\BC{\begin{corollary}}
\def\EC{\end{corollary}}
\def\BE{\begin{example}}
\def\EE{\end{example}}
\def\BD{\begin{definition}}
\def\ED{\end{definition}}
\def\BR{\begin{remark}}
\def\ER{\end{remark}}
\def\BAS{\begin{assumption}}
\def\EAS{\end{assumption}}
\def\BI{\begin{itemize}}
\def\EI{\end{itemize}}

\def\la{\langle}
\def\ra{\rangle}

\def\argmin{\mathop{\rm argmin}}

\newcommand{\al}{\alpha}

\newcommand{\g} {\gamma}
\newcommand{\wh}[1] {\widehat{#1}}

\newcommand{\dl}{\delta}

\def\bE{{\mathbf{E}}}

\newcommand{\e}{\epsilon}

\newcommand{\rl}[1]{~\cite{#1}}

\newcommand{\sign}{\mbox{sign}}

\newcommand{\cF}{{\cal F}}

\newcommand{\cA}{{\cal A}}
\newcommand{\cC}{{\cal C}}

\newcommand{\cB}{{\cal B}}

\newcommand{\be}{\begin{eqnarray}}
\newcommand{\ee}[1]{\label{eq:#1}\end{eqnarray}}
\newcommand{\nn}{\nonumber \\}
\newcommand{\ese}{\end{eqnarray*}}
\newcommand{\bse}{\begin{eqnarray*}}
\newcommand{\rf}[1]{~(\ref{eq:#1})}

\newtheorem{theo}{Theorem}
\newcommand{\bthm}{\begin{theo}}
\newcommand{\ethm}[1]{\label{the:#1}\end{theo}\par}
\newcommand{\rth}[1]{\ref{the:#1}}
\newtheorem{col}{Corollary}
\newcommand{\bcol}{\begin{col}}
\newcommand{\ecol}[1]{\label{the:#1}\end{col}\par}
\newcommand{\rc}[1]{\ref{the:#1}}
\newtheorem{defi}{Definition}
\newcommand{\bdf}{\begin{defi}}
\newcommand{\edf}[1]{\label{df:#1}\end{defi}\par}

\newtheorem{lem}{Lemma}
\newcommand{\blem}{\begin{lem}}
\newcommand{\elem}[1]{\label{le:#1}\end{lem}\par}

\newtheorem{pro}{Proposition}
\newcommand{\bpro}{\begin{pro}}
\newcommand{\epro}[1]{\label{pr:#1}\end{pro}\par}
\newcommand{\rp}[1]{\ref{pr:#1}}
\newtheorem{statement}{Problem}
\newcommand{\bstatement}{\begin{statement}}
\newcommand{\estatement}[1]{\label{stat:#1}\end{statement}\par}

\newcounter{assc}
\newcommand{\bass}[1]{\refstepcounter{assc}\label{ass:#1} \begin{maliste}{\bf \arabic{assc}.}}
\newcommand{\eass}{\end{maliste}}
\newcommand{\ras}[1]{\ref{ass:#1}}
\newcounter{algoc}
\newcommand{\balgo}[1]{\refstepcounter{algoc}\label{algo:#1} \begin{malistea}{\bf \arabic{algoc}.}}
\newcommand{\ealgo}{\end{malistea}}
\newcommand{\ralgo}[1]{\ref{algo:#1}}
\newcommand{\pr} {\noindent{\bf Proof\,:~}}

\def    \qed    {\hfill\hbox{\hskip 4pt
                \vrule width 3pt height 1pt depth 1pt
                \hbox{\vrule width 2pt height 3.5pt depth 1pt}}}

\newcommand{\aic}[2]{{\color{blue} #1}{\color{red} ~#2}}
\begin{document}
\title{Primal-dual subgradient methods
for minimizing uniformly convex functions}
\author{Anatoli Juditsky\thanks{LJK,
Universit\'e J. Fourier, B.P. 53, 38041 Grenoble
Cedex 9, France, {\tt juditsky@imag.fr}},\and
{Yuri} {Nesterov}\thanks{CORE,
Catholic University of Louvain, 34 voie du Roman Pays,
 1348 Louvain-la-Neuve, Belgium,
{\tt nesterov@core.ucl.ac.be}}}
\maketitle
\begin{abstract}
We discuss non-Euclidean deterministic and stochastic algorithms for optimization problems with strongly and uniformly convex objectives. We provide accuracy bounds for the performance of these algorithms and design methods which are adaptive with respect to the parameters  of strong or uniform convexity of the objective: in the case when the total number of iterations $N$ is fixed, their accuracy coincides, up to a logarithmic in $N$ factor with the accuracy of optimal algorithms.
\end{abstract}

\section{Introduction}

\vspace{1ex}\noindent
Let $E$ be a (primal) finite-dimensional real vector space.
In this paper we consider the optimization
problem:
\beq\label{prob-Main}
\min\limits_{x} \{ f(x):\; x \in Q \},
\eeq
where $Q$ is a closed convex set in $E$ and function $f$
is {\em uniformly convex} and Lipschitz-continuous on $Q$. Recall that
a function $f$ is called {\em uniformly convex} on $Q
\subset E$ with convexity parameters $\rho=\rho(f)\ge 2$ and  $\mu = \mu(f,\rho)$  if for
all $x$ and $y$ from $Q$ and any $\alpha \in [0,1]$ we
have
\be
f(\alpha x + (1-\alpha)y) &\leq& \alpha f(x) + (1-\alpha)
f(y) \nn
&&- \half \mu \alpha(1-\alpha) [\alpha^{\rho-1}+(1-\alpha)^{\rho-1}]\| x - y \|^\rho.
\ee{def-Conv0}
The function $f$ which is uniformly convex with $\rho=2$ is called {\em strongly convex}.
Uniform convexity with $2\le \rho \le \infty$ and $\mu\ge 0$ implies usual convexity.
\par
In this paper we discuss deterministic and stochastic first order algorithms for (large scale) {\em non-Euclidean uniformly convex objectives}, thus extending  non-Euclidean first order methods  (see, e.g. \cite{NJLS,DA} and references therein) to  uniformly convex optimization.
\par
Uniformly convex functions have been introduced to optimization in \rl{Polyak} and extensively studied (cf. \rl{AP},
\rl{VNC1}, and \rl{Zali}).  The worst-case complexity bounds for the problem (\ref{prob-Main}) with the exact and stochastic  first order oracle are available for the case of strongly convex objective (see, e.g. \cite{RR2011,ABRW2102} and references therein). Specifically, for any method tuned to the absolute accuracy $\e$ for the problem (\ref{prob-Main}) with strongly convex, with parameter $\mu$, and Lipschitz-continuous, with unit Lipschitz constant, objective and deterministic first order oracle, the number of calls to the oracle is not less than $O(\mu^{-1}\e^{-1})$ which is much better than the corresponding bound $O(\e^{-2})$ for a  larger class of Lipschitz-continuous convex functions.
The corresponding bound for uniformly convex problems with the convexity parameters $\rho$ and $\mu$ reads $O\left(\mu^{-{2\over \rho}} \e^{-{2(\rho-1)
\over \rho}}\right)$ (for the sake of completeness we provide in appendix \ref{app-lower} the corresponding bound for the case of  the Euclidean norm  $\|\cdot\|$). Note that in the case of the stochastic oracle these bounds holds also for problems with smooth objective.
\par
Note that {\em smooth uniformly convex deterministic} optimization is ``covered'' within the Euclidean framework -- it appears that  the optimal deterministic first order algorithms of Euclidean {\em smooth uniformly convex optimization} developed in \cite[chapter 7]{nemyud:83} and \cite[chapter 2]{Nesterov-book} retain their optimality in the non-Euclidean framework. Indeed, let us consider the problem (\ref{prob-Main}) where $f$ is a strongly convex quadratic form:
$
f(x)=\half x^T Ax -b^Tx$,
the set $Q=\{x\in \mathbb R^n|\,\|x\|_1\le 1\}$, and $A$ is a  symmetric $n\times n$ positive-definite matrix. Recall that the complexity estimate for optimal algorithms of strongly convex smooth optimization is $O(\sqrt{\lambda}\log \epsilon^{-1} )$ where
$\lambda={{\cal L}(f)\over \mu(f)}$ is the conditioning of the objective -- the ratio of the Lipschitz constant ${\cal L}(f)$ of the gradient of the objective and the parameter $\mu(f)$ of strong convexity, and $\epsilon$ is the desired absolute accuracy.
Note that the Lipschitz constant ${\cal L}_1(f)$ of the gradient of $f$ with respect to the norm $\|\cdot\|_1$ satisfies
${\cal L}_1(f)=\|A\|_{1,\infty}=\max_{1\le i,j\le n}|A_{ij}|$. On the other hand, one may easily verify that the corresponding parameter $\mu(f)$ of strong convexity of $f$ is bounded from above with ${\cal L}_1(f)n^{-1}$, resulting in conditional number $\lambda\ge n$.\footnote{Here is the proof of this claim: let $\xi=(\xi_1,...,\xi_n)^T$ be a random vector with i.i.d. components such that $P(\xi_i=1/n)=P(\xi_i=-1/n)=1/2$. Then $\|\xi\|_1=1$, and
\[\mu(f)\le E(\xi^TA\xi)={1\over n^2} \sum_{i}A_{ii}\le  {{\cal L}_1(f)\over n}.
 \]Observe that the bound ${1/n}$ is attained for the identity matrix $A$.
}
 Now recall that the Lipschitz constant of the gradient of $f$, when measured  with respect to Euclidean norm is
${\cal L}_2(f)=\|A\|_{2,2}=\lambda_{\max}(A)$ -- the spectral norm of $A$, and ${\cal L}_2(f)\le n{\cal L}_1(f)$. In other words, in this case, when passing from the Euclidean to non-Euclidean setup we gain nothing -- the degradation of the strong convexity parameter in the $\|\cdot\|_1$-setup outweighs the potential improvement of the conditioning due to the reduced Lipschitz constant in the $\|\cdot\|_1$-setup.
\par
On the other hand, although the optimal algorithms for optimization with {\em strongly convex Lipschitz continuous} objective in the Euclidean framework are readily available (see, e.g., \cite{NV,RR2011}), they cannot be directly transposed to the non-Euclidean framework.
\par
The results presented in this paper are not very new, as they were developed by the authors in 2004-2005. However, because of the immediate lack of application and, more importantly, due to new first order methods based on smoothing of structured problems with better complexity characteristics
which were developed in \cite{SMNF,EG} at that time, the authors got an impression that new non-Euclidean algorithms of black-box (non-structured)
uniformly convex optimization are of very limited interest.
However, certain developments of the last years clearly
demonstrated that in some situations the black-box methods are irreplaceable. Indeed,
 exact first order oracle are often unavailable, or the structure of a problem may be simply too complex for applying a smoothing technique.
In particular, deterministic and stochastic non-Euclidean first order methods of convex optimization have attracted much attention lately in relation, in particular, with very large scale applications arising in statistics and learning. For instance, some new applications involving large scale strongly convex optimization has been recently reported (see, e.g., \cite{Srebro,Xiao,Pag}). These considerations encouraged the authors to publish the above mentioned results
on subgradient methods for uniformly convex problems.
\par
In this paper we
develop minimax optimal primal-dual minimization schemes in the spirit of
\cite{DA} for uniformly convex problems as in (\ref{prob-Main}) with Lipschitz-continuous objective.
We also study the performance of multistage dual averaging procedures when applied to uniformly convex stochastic
minimization problems. In particular, we show that such procedures attain the minimax rates of convergence on the considered problem class. We also provide confidence sets for approximate solutions of stochastic uniformly convex problems.
\par
It is well known that performance of ``classical'' optimization routines for strongly (and uniformly) convex problems can become very poor when the parameters of strong (uniform) convexity are not known {\em a priori} (see, e.g. section 2.1 in \cite{NJLS}). In the case of deterministic and stochastic optimization we develop {\em adaptive} minimization procedures in the case when the total number $N$ of the method iterations is fixed. The accuracy of these procedures (which do not require a priori knowledge of parameters of uniform convexity) coincides, up to a logarithmic in $N$ factor, with the accuracy of optimal algorithms (which ``know'' the exact parameters). It is worth to note that we do not know if it is possible to construct adaptive optimization procedures tuned to the fixed accuracy with analogous proprieties.
\par
The paper is organized as follows: in section \ref{sec-probst} we define the basic ingredients of the minimization problem in question. Then we study the properties of the primal-dual subgradient algorithms in the problem with an exact deterministic oracle in section \ref{deterministic} and show how the dual solutions can be produced in section \ref{sc-PD}. In section \ref{stoch} we develop optimal algorithms  for stochastic uniformly convex optimization and show how confidence sets for approximate solutions can be constructed. Section \ref{sec-comp} contains some details of computation aspects of proposed routines. Finally, in appendix \ref{app-lower} we present the lower complexity bound for a class of optimization problems with uniformly convex and Lipschitz continuous objectives; appendix \ref{app-proofs} contains the proofs of the statements of the paper.
\section{Problem statement and basic assumptions}
\label{sec-probst}
\subsection{Notations and generalities} Let $E^*$ be t
he dual of $E$. We denote the value of linear function $s \in E^*$
at $x \in E$ by $\la s, x \ra$. For measuring distances in
$E$, let us fix some (primal) norm $\| \cdot \|$. This
norm defines a primal unit ball
$$
B = \{ x \in E: \; \| x \| \leq 1 \}.
$$
The dual norm $\| \cdot \|_*$ on $E^*$ is introduced, as
usual, by
$$
\| s \|_* = \max\limits_x \{ \la s, x \ra: \; x \in B \},
\quad s \in E^*.
$$
For other balls in $E$ we adopt the following notation:
$$
B_R(x) = \{ y \in E:\; \| y - x \| \leq R \}, \quad x \in
E.
$$
%
If a uniformly convex function $f$ is {\em subdifferentiable} at $x$, then
\beq\label{def-Conv1}
f(y) \geq f(x) + \la f'(x), y - x \ra + \half \mu \| y - x
\|^\rho \quad \forall y \in Q,
\eeq
where $f'(x) \in E^*$ denotes one of {\em subgradients} of
$f$ at $x \in Q$. If $f$ is subdifferentiable at two
points $x, y \in Q$, then\footnote{\aic{Note that the relationship (\ref{def-Conv12}) is sometimes used as definition of a uniformly convex function (see, e.g. \cite{Nest:2008}). However, (\ref{def-Conv12}) does not imply (\ref{def-Conv1}) and  \rf{def-Conv0}, but, instead of (\ref{def-Conv1}), for instance, it leads to
\[
f(y) \geq f(x) + \la f'(x), y - x \ra + {\mu\over \rho} \| y - x
\|^\rho \quad \forall y \in Q.\]
Of course, in the strongly convex case we have $\rho=2$ and both definitions lead to the same value of the modulus of strong convexity. }{}}
\beq\label{def-Conv12}
\la f'(x) - f'(y), x - y \ra \geq \mu \| x - y \|^\rho.
\eeq

%
%
\subsection{Problem statement}
\label{prstat}
We consider the optimization problem (\ref{prob-Main}) with the uniformly convex function $f$
with convexity parameters $\rho(f)$ and $\mu(f)$. The basic assumption we make about the objective,
and which is supposed to hold through the paper, is that $f$ is Lipschitz-continuous on $Q$:
\bass{det1}
We assume that all subgradients of the objective function are bounded:
\[
\|f'(x)\|_*\le L, \;\;\;\mbox{for any}\;\; x\in Q.
\]
\eass
We are to study the performance of an iterative minimization schemes, and we consider
two settings which differ with respect to the information available to the method at each iteration.
\begin{description}
\item[] {\em -- deterministic setting:} let $x_k$ be the search points at iteration $k$, $k=0,1,...$.
We suppose that
an exact subgradient observations $g_k=f'(x_k)$ and the exact objective values $f(x_k)$ are available;
\item[] {\em -- stochastic setting:} the observation $g_k$ of the subgradient $f'(x_k)$, requested by the method
at the $k$-th iteration,
is supplied by a {\sl stochastic oracle}, i.e. $g_k$ is a random vector.
\end{description}
To be more precise,
suppose that we are given the probability space $(\Omega, \cF, P)$ and a filtration
$(\cF_k)$, $k=-1,0,1,...$
(non-decreasing family of $\sigma$-algebras which satisfies ``usual" conditions).
\par
Let
\[
g_k\equiv g(x_{k},\omega_k),
\]
where
\begin{itemize}
     \item $\{\omega_k\}_{k=0}^\infty$ is sequence of random parameters taking
values in $\Omega$, such that $\omega_k$ is $\cF_{k}$-measurable;
     \item  $x_k$ is the $k$-th search point generated by the method. We suppose that
$x_k$ is $\cF_{k-1}$-measurable (indeed, $x_k$ is a measurable function of
$x_0$ and observations $g_1,...,g_{k-1}$ at iterations
$1,...,k-1$).
\end{itemize}
We also consider the following assumptions specific to the stochastic problem:
\bass{stoch1}
The oracle is unbiased. Namely,
\bse
\bE_{k-1}[g(x_k,\omega_k)]\in
\partial f(x_k),\;\;\mbox{a.s.}\;\;\; x_k\in Q,\;\;k=0,1,...
\ese
\eass
Here $\bE_{k}$ stands for the expectation conditioned by $\cF_{k}$  (then $\bE=\bE_{-1}$ is the ``full'' expectation).

Let us denote
\[
\xi_k=g_k-f'(x_k),
\]
the stochastic perturbation. Note that $\bE_{k-1}[\xi_k]=0$ a.s. for $k=0,1,...$.
We suppose that the intensity of the sequence $\{g_k\}_{k=0}^\infty$ is bounded.
\bass{stoch2}
We assume that
\be
\sup_k \bE\|\xi_k\|^2_*\le \sigma^2<\infty\;\;\mbox{for }\;\; k=0,1,...
\ee{sigma2}
\eass

 We will also use a stronger bound on the tails of the distribution of $(\xi_k)$:
\bass{stoch3}
There exists $\sigma<\infty$ such that
\be
\bE_{k-1}\left[\exp\left\{{\|\xi_k\|_*^2\sigma^{-2}}\right\}\right]\le \exp(1)\;\;\mbox{a.s.}, \;\;k=0,1,...
\ee{expsigma}
\eass
Note that by the Jensen inequality \rf{expsigma} implies \rf{sigma2}.
\subsection{Prox-function of the unit ball
}\label{sc-Ball}

Assume that we know a {\em prox-function} $d(x)$ of the
ball $B$. This means that $d$ is continuous and strongly
convex on $B$ in terms of \rf{def-Conv0} with some
convexity parameter $\mu(d) > 0$. Moreover, we assume that
$$
d(x) \geq d(0) = 0, \quad x \in B.
$$
Hence, in view of (\ref{def-Conv1}) we have
\[
d(x) \geq \half \mu(d) \| x \|^2, \quad \forall x \in
Q\cap B.
\]
An important characteristic of the prox-function is
its maximal value on the unit ball:
\be
d(x) \; \leq \; A(d), \quad x \in B.
\ee{eq-UpBound}
Therefore,
\beq\label{eq-AMu}
\ba{rcl}
\mu(d) & \leq & 2 A(d).
\ea
\eeq
If the function $d$ is growing quadratically, another important characteristics
is its constant of quadratic growth $C(d)$ which we define as the smallest $C$ such that
\be
d(x)\le C\|x\|^2.
\ee{quadgr}
We have
\[
\mu(d)\le 2C(d)\;\;\;\mbox{and}\;\;\; A(d)\le C(d).
\]
\paragraph{Example 1.}  Let $E=\mathbb R^n$ and let $B$ be a unit Euclidean ball in $\mathbb R^n$.
We choose the norm $\|\cdot\|$ to be the Euclidean norm on $\mathbb R^n$, so that the function
$d(x)=\|x\|_2^2/2$ is strongly convex with $\mu(d)=1$ and $C(d)=A(d)=1/2$.
\paragraph{Example 2.} Let again  $E=\mathbb R^n$ and let $B$ be the standard hyperoctahedron in $\mathbb R^n$, i.e. a unit $l_1$-ball: $B=\{x\in \mathbb R^n| \,\|x\|_1\le 1\}$, where
\[
\| x \|_1 = \sum\limits_{i=1}^n | x^{(i)} |.
\]
We take $\|x\|=\|x\|_1$ and consider for $p>1$ the function $d$,
\[
d(x)=\half \left(\sum_{i=1}^n |x_i|^p
\right)^{2/p}=\half {\|x\|_p^2}.
\]
The function $d$ is strongly convex with  $\mu(d)=O(1)n^{{p-1\over p}}$, and for $p=1+{1\over \ln n}$ we have $\mu(d)=O(1)(\ln n)^{-1}$ (see, e.g. \cite{nemyud:83}). Further, we clearly have $A(d)=C(d)=1/2$.
\par
Note that norm-type prox-functions are not the only possible in the hyperoctahedron setting. Another example of prox-function of the $l_1$-unit ball $B$, which is very interesting from the computational point of view, is as follows:
%
\beq\label{eq:DistL1}
\ba{rl}
d(x) \; = & \min \{\;  \sum\limits_{i=1}^n [ \;
\psi(u^{(i)})  +  \psi(v^{(i)}) \; ]: \;
\sum\limits_{i=1}^n \left[ u^{(i)}+ v^{(i)} \right] = 1,\\
\\ &  x^{(i)} = u^{(i)} - v^{(i)},\; u^{(i)} \geq 0,\;
v^{(i)} \geq 0,  \; i = 1, \dots, n \; \} \; + \;
\ln (2n), \\
\\
\psi(t ) \; = & \left\{ \ba{rl} t \ln t, & t > 0, \\ 0, &
t = 0 . \ea \right.
\ea
\eeq
In order to show that this function is strongly convex
on the standard hyperoctahedron $B=\{x\in \mathbb R^n|\,
\|x\|_1\le 1\}$, we need the following general result.
\BL\label{lm-Symm} Let $Q$ be a bounded closed convex set in $E$
containing the origin. If function $f(x)$ is strongly
convex on $Q$ with parameter $\mu \geq 0$, then its
symmetrization
$$
\ba{rcl}
f^0(x) & = & \min\limits_{u,v, \alpha} \left\{ f(u) +
f(v):\; x = u - v, \; u \in \alpha Q, \; v \in (1 -
\alpha)Q,\; \alpha \in [0,1] \right\},
\ea
$$
is strongly convex on the set $Q^0 = \conv \{ Q, - Q \}$
with convexity parameter $\half \mu(f)$.
\EL
Thus, for function $d(x)$ defined by \rf{DistL1} we
can take
\[
\mu(d) = \half, \quad A(d) = \ln (2n).
\]
Note that $d$ does not satisfy the quadratic growth condition \rf{quadgr}.
\par
For $z\in Q$, consider the set
\bse
Q_{R}(z) \Def Q\cap B_R(z).
\ese
This set can be equipped with a prox-function
\bse
d_{z, R}(x) = d\left({1 \over R}(x - z) \right).
\ese
Thus, the prox-center of the set $Q_{R}(z)$ is $z$, and
$\mu(d_{z, R}) = {1 \over R^2} \mu(d)$. Moreover, by
\rf{eq-UpBound},
\bse
d_{z, R}(x) \leq A(d), \quad \forall x \in Q_{R}(z).
\ese

In what follows we need the objects: the function
\be
V_{z, R, \beta}(s) & = & \max\limits_{x} \{ \la s, x -
z\ra - \beta d_{z, R}(x):\; x \in Q_{R}(z) \},
\ee{V}
and the prox-mapping
\bse
\pi_{z, R, \beta}(s) & = & \arg\max\limits_{x} \{ \la s, x
- z \ra - \beta d_{z, R}(x):\; x \in Q_{R}(z) \}.
\ese
Note that $\dom V_{z,R,\beta}=E^*$. Let us mention some
properties of function $V_{z,R,\beta}$ (cf. Lemma 1
\rl{DA}):
\begin{itemize}
\item if $\beta_1\le \beta_2$ then $V_{z, R,\beta_1}(s)\ge V_{z, R,\beta_2}(s)$;
\item the function $V_{z,R,\beta}$ is convex and differentiable on $E^*$. Moreover,
its gradient is Lipschitz continuous with the constant
${R^2 \over \beta \mu(d)}$:
\bse
\ba{rcl}
\|V'_{z,R,\beta}(s_1)-V'_{z,R,\beta}(s_2)\| & \le &
{R^2\over
\beta \mu(d)} \|s_1-s_2\|_*,\quad \forall s_1, s_2
\in E^*.
\ea
\ese
\item
For any $s\in E^*$,
\bse
V'_{z,R, \beta}(s)+z=\pi_{z,R,\beta}(s)\in Q_{R}(z).
\ese
\end{itemize}
\section{Deterministic methods for uniformly convex functions}
\label{deterministic}
We start with the description of the basic tool -- the dual averaging procedure, which originates
in \rl{DA}.
\subsection{Method of Dual Averaging}
%
At each phase the dual averaging (DA) method will be applied to the following auxiliary
problem: 
\beq\label{prob-Aux}
\min\limits_{x} \{ f(x):\; x \in Q_{R}(\bar x) \}.
\eeq
Its feasible set is endowed with the following
prox-function:
\[
\ba{c}
d_{\bar x, R}(x) = d\left({1 \over R}(x - \bar x) \right).
\ea
\]
Consider now the generic scheme of Dual Averaging as
applied to the problem (\ref{prob-Aux}).
\balgo{DA-met}
\begin{description}
\item[Initialization:] Set $x_0 = \bar x$, $s_0 =
0 \in E^*$.
Choose $\beta_0 > 0$.
\item[Iteration] ($k \geq 0$):
\begin{enumerate}
\item Choose $\lambda_{k} > 0$. Set $s_{k+1} =
s_{k}+ \lambda_{k} f'(x_{k})$, where $\{ \lambda_i \}_{i=0}^{\infty}$ is a sequence of
positive parameters.
\item Choose $\beta_{k+1} \geq \beta_{k}$. Set
$x_{k+1} = \pi_{\bar x, R, \beta_{k+1}}(- s_{k+1})$.
\end{enumerate}
\end{description}
\ealgo
The process is terminated after $N$ iterations. The
resulting point is defined as follows:
\beq\label{eq-AuxS}
\ba{c}
x_N(\bar x, R) \; = \; \left(\sum_{i=0}^N\lambda_i\right)^{-1} \sum\limits_{i=0}^N\lambda_i
x_i.
\ea
\eeq
The result below underlies  the following developments (cf. Theorem 1 of \cite{DA}.):
\bpro
 For any $x\in Q_R(\bar{x})$,
\be
\sum_{i=0}^k \lambda_i\la f'(x_i),x_i-x\ra \leq
\;d_{\bar{x},R}(x)
\beta_{k+1}  +{R^2
\over 2 \mu(d)} \sum\limits_{i=0}^k {\lambda_i^2 \over
\beta_i} \|f'(x_i) \|_*^2.
\ee{eq-DADelta}
\epro{DA}
\par
Let $\lambda_i=1$ and
$\beta_i = \gamma \sqrt{N+1}$, $i=0,...,N$ with some $\gamma > 0$.
%
We form the {\em gap value}

\be
\delta_k(\bar x, R) = \max\limits_{x} \left\{{1\over k+1}
\sum\limits_{i=0}^k \la f'(x_i), x_i - x \ra:\;
x \in Q_{R}(\bar x) \right\}.
\ee{gapv}
In view of (\ref{eq-AuxS}) we have the
following lemma:
\BL\label{lm-DA}
Let us choose an arbitrary $\bar x \in Q$
and let $x^*$ be the optimal
solution of problem (\ref{prob-Aux}).
Then the
approximate solution supplied by Algorithm \ralgo{DA-met} with the constant gain
$\beta_i=\gamma\sqrt{N+1}$
satisfies
\bse
f(x_N(\bar x,R)) - f(x^*) & \leq &  {1 \over \sqrt{N+1}}
\left( \gamma A(d) + {L^2 R^2 \over 2 \gamma \mu(d)}
\right),\\
\| x_N(\bar x,R) - x^* \|^\rho & \leq & {\dl_N(\bar{x},R)\over \mu(f)}\le {1 \over \mu(f)
\sqrt{N+1}} \left( \g A(d) + {L^2 R^2 \over 2 \gamma
\mu(d)} \right).
\ese
\EL
%
Under the premises of the lemma we can establish the following immediate bounds:
\BC\label{the:cor-Stage1}
Let $x^*$ be an optimal solution of (\ref{prob-Aux}). Then for the choice
\[
\gamma = { L R \over \sqrt{2 \mu(d) A(d)}}
\]
we have the estimates:
\begin{eqnarray}\label{eq:eq-FB1}
f(x_N(\bar x, R)) - f(x^*) \; &\leq &\;L R \sqrt{2 A(d) \over \mu(d) (N+1)},
\\
\| x_N(\bar x, R) - x^* \|^\rho & \leq & {L R \over \mu(f)}\sqrt{2 A(d)
\over\mu(d) (N+1)}.\nonumber
\end{eqnarray}
\EC

\subsection{Multi-step algorithms }
\label{sec:smallballs}
Now we are ready to analyze multistage procedures for uniformly convex functions. In this section we assume
that the constants $L$, $\mu(f)$, $\rho$ and $R_0 \geq \| x^* -
x_0 \|$ are known. Let us fix $\e>0$ and let  $x_0$ be an arbitrary element of $Q$.
\balgo{1}
\begin{description}
\item[Initialization:]
Set $y_0 = x_0$ and $m =\lfloor\log_2 {\mu(f) \over
\epsilon} R_0^\rho\rfloor+1$.\footnotemark
\footnotetext[1]{$^)$Here $\lfloor a \rfloor$ stands for the
largest integer strictly smaller than $a$.}$^)$ Let $\tau={2(\rho-1)\over \rho}$.
\item[Stage] $k = 1, \dots , m$:
\begin{enumerate}
\item Define $N_k = \lfloor 2^{\tau k} {4 L^2 A(d)\over \mu^2(f)
\mu(d) R_0^{2(\rho-1)}}\rfloor$ and $R^\rho_k = 2^{-k}R^\rho_0$.
\item Compute $y_k = x_{N_k}(y_{k-1},R_{k-1})$ with
$\gamma_k = { LR_{k-1} \over \sqrt{ 2 \mu(d) A(d) }}$.
\end{enumerate}
\item[Output:] $\wh{x}_\e(y_0,R_0):=y_m$.
\end{description}
\ealgo

Note that the parameters of the algorithm satisfy the
following relations:
\begin{equation}\label{eq:eq-AuxD}
N_k + 1 \; \geq \;  2^{\tau k} {4 L^2 A(d)\over \mu^2(f) \mu(d)
R_0^{2(\rho-1)}} \; \geq \; N_k , \quad\quad 2^m \; \geq \; {\mu(f)
\over \epsilon} R_0^\rho \; \geq \; 2^{m-1}.
\end{equation}
\bthm
The points $\{ y_k \}_{k=1}^m$ generated by
Algorithm \ralgo{1} satisfy the following conditions:
\beq\label{eq-Rad}
\ba{rcl}
\| y_k - x^* \|^\rho & \leq & R^\rho_k = 2^{-k}R^\rho_0, \quad k = 0,
\dots, m,
\ea
\eeq
\beq\label{eq-Val}
\ba{rcl}
\dl_{N_{k}}(y_{k-1},R_{k-1})& \le & \mu(f) R_k^\rho \; = \;
\mu(f) 2^{-k} R_0^\rho, \quad k = 1, \dots, m.
\ea
\eeq
Moreover, $f(\wh{x}_{\epsilon}(y_0,R_0)) - f^* \leq
\epsilon$ and the total number $N(\e)$ of iterations in the
scheme does not exceed
\beq\label{eq-IT}
\ba{rcl}
\left({2^{{m+1}}\over R_0^\rho}\right)^\tau {4 L^2 A(d)\over \mu^2(f) \mu(d)}
&\stackrel{\rf{eq-AuxD}}{\leq} & {4^{\tau+1} L^2 A(d)\over
\mu(f)^{2\over \rho} \mu(d)} \, \epsilon^{-\tau}.
\ea
\eeq
\ethm{th-Mult1}
An important particular case of Theorem \rth{th-Mult1} is the case of strongly convex objective $f$.
In the latter case $\tau=1$ and the analytical complexity of Algorithm \ralgo{1}
does not exceed
\[
{16 L^2A(d)\over \mu(f)\mu(d)}\,\e^{-1}.
\]
The method can be easily rewritten for the case when the total
number $N$ of calls to the oracle is fixed a priori.\par
Denote $\bar{N}=[2^{\tau}(2^\tau+1)] {4L^2A(d)\over \mu^2(f)\mu(d) R_0^{2(\rho-1)}}$. If $N<\bar{N}$
run Algorithm 1 with $
\gamma = { L R_0 \over \sqrt{2 \mu(d) A(d)}}$ and output the approximate solution $\wh{x}=x_N(\bar{x},R_0)$.
If $N\ge \bar{N}$ use the following procedure:
\balgo{2}
\begin{description}
\item[Initialization:] set $y_0 = x_0$, $\tau={2(\rho-1)\over \rho}$, compute
$N_j=\lfloor 2^{\tau j}{4L^2A(d)\over \mu(f)^2\mu(d)R_0^{2(\rho-1)}}
\rfloor$ while $\sum_{j} N_j\le N$. Set
$$
\ba{rcl}
m(N)& = & \max\{k:\;\sum_{j=1}^k N_j\le N\}.
\ea
$$
\item[Stage] $k = 1, \dots , m(N)$:
Set $R^\rho_k=2^{-k}R^\rho_0$. Compute $y_k =
x_{N_k}(y_{k-1},R_{k-1})$ with
$$
\ba{rcl}
\gamma_k & = & {LR_{k-1} \over \sqrt{ 2 \mu(d) A(d) }}.
\ea
$$
\item[Output:] $\wh{x}_N=y_{m(N)}$.
\end{description}
\ealgo

\BC\label{mult2}
We have
\beq\label{eq-it2}
\ba{rcl}
f(\wh{x}_N) - f^* & \leq & 2\left({8L^2A(d)\over \mu(f)^{2\over \rho}\mu(d)N}\right)^{1/\tau}.
\ea
\eeq
\EC

\subsection{Methods with quadratically growing prox-function}
We propose here a slightly different version of multi-stage procedures for the case
when the prox-function satisfies the condition \rf{quadgr} of quadratic growth.

\par
The result below is an immediate consequence of Proposition \rp{DA} (cf. Lemma \ref{lm-DA}
and Corollary \rc{cor-Stage1}):
\BC\label{the:cor-Stage2}
Let $x^*$ be an optimal solution of (\ref{prob-Aux}).
Suppose that the prox-function $d$ satisfies \rf{quadgr} and that
$\|\bar{x}-x^*\|\le r\le R$.
Then the approximate solution $x_N(\bar x, R)$, provided by Algorithm \ralgo{DA-met} with
\[
\gamma={R^2L\over r\sqrt{2C(d)\mu(d)}},
\]
satisfies
\be
\label{eq:eq-FB0}
f(x_N(\bar x, R)) - f(x^*) \; &\leq& \; rL\sqrt{2C(d)\over \mu(d)(N+1)},\\
\| x_N(\bar x, R) - x^* \|^\rho & \leq & {r L \over \mu(f)}\sqrt{2 C(d)\over \mu(d) (N+1)}.
\ee{eq-NB0}
\EC
Indeed, to show  \rf{eq-FB0} and
\rf{eq-NB0} it suffices to use \rf{eq-DADelta} and to observe that due to \rf{quadgr}
$d_{\bar x,R}(x^*)\le C(d){r^2\over R^2}$.
\par
The following multi-stage scheme exploits the ``scalability property'' \rf{quadgr}
of the prox-function $d$.
It starts from arbitrary $x_0 \in Q$. As in the previous section, we assume
that the constants $L$, $\mu(f)$ and the diameter $R_0$ of
$Q$ are known.
\balgo{5}
\begin{description}
\item[Initialization:]
Set $y_0 = x_0$, $\tau={2(\rho-1)\over \rho}$ and $m =\lfloor\log_2 {\mu(f) \over
\epsilon} R_0^\rho\rfloor+1$.
\item[Stage] $k = 1, \dots , m$:
\begin{enumerate}
\item Define $N_k = \lfloor 2^{\tau k} {4 L^2 C(d)
\over \mu^2(f)
\mu(d) R_0^{2(\rho-1)}}\rfloor$ and $r^\rho_k = 2^{-k}R^\rho_0$.
\item Compute $y_k = x_{N_k}(y_{k-1},R_0)$ with
$\gamma_k = {LR_0^2\over r_{k-1}
\sqrt{2C(d)\mu(d)}}$.
\end{enumerate}
\item[Output:] Set the approximate solution $\wh{x}_\e=y_m$.
\end{description}
\ealgo
We would like to stress the difference between
Algorithms \ralgo{1} and \ralgo{5}: in
Algorithm \ralgo{5} the delation parameter $R=R_0$ of the
prox-function $d$ remains the same through all the stages
of the method. Only the gain  $\gamma_k$ and the duration
$N_k$ of the stage depend on the stage index $k$. As a result, the prox-mapping $ \pi_{z, R, \beta}$ is easier to compute.
Further, as we will see in section \ref{sec:expect}, it also allows a straightforward modification in the case of stochastic oracle.
\par
We have
the following  analogue of Theorem
\rth{th-Mult1} in this case:
\bthm Suppose that
\[
N\ge N(\e)={4^{\tau+1} L^2 C(d)\over \mu(f)^{2\over \rho}\mu(d)}
\e^{-\tau}.
\]
Then the approximate solution $\wh{x}_N$, provided by
Algorithm \ralgo{5} satisfies:
\[
f(\wh{x}_\e)-f^*\le \e.
\]
\ethm{multi-qdet}
The method can be rewritten when
the total number $N$ of calls to the oracle is fixed.
\par
Suppose that $N\ge 2^\tau(2^\tau+1) {4L^2C(d)\over
\mu^2(f)\mu(d) R_0^{2(\rho-1)}}$. Consider the following procedure:
\balgo{22}
\begin{description}
\item[Initialization:] Set $y_0 = x_0$, $\tau={2(\rho-1)\over \rho}$, compute
$N_j=\left\lfloor 2^{\tau j}{4L^2C(d)\over
\mu(f)^2\mu(d)R_0^{2(\rho-1)}}\right\rfloor$, while $\sum_{j} N_j\le
N$. Set $m(N)=\max\{k:\;\sum_{j=1}^k N_j\le N\}$.
\item[Stage] $k = 1, \dots , m(N)$:\\
Set $r_k^\rho=2^{-k}R^\rho_0$. Compute $y_k =
x_{N_k}(y_{k-1},R_0)$ with $\gamma_k =
{LR_0^2\over r_{k-1} \sqrt{2C(d)\mu(d)}}$.
\item[Termination:] Set the approximate solution $\wh{x}_N=y_{m(N)}$.
\end{description}
\ealgo
\BC \label{the:mult2-det}
We have
\[
 f(\wh{x}_N) - f^* \leq 2\left({8L^2C(d)\over \mu(f)^{2\over \rho}\mu(d)N}\right)^{1/\tau}.
\]
\EC
The proof of the corollary is completely analogous to that of Corollary \ref{mult2}.
%
\subsection{Adaptive algorithm}
\label{adaptive}
Consider the setting in which the total number $N$
of calls to the oracle is fixed and suppose that the convexity
parameters $ \rho$ and $\mu(f)$ are unknown.
We propose a multi-stage procedure
which does not require the knowledge of these parameters and attains the accuracy of
the method which ``knows" the convexity parameters up
to a logarithmic in $N$ factor. Following the terminology used in statistics and control literature,
we call such procedures adaptive (with respect to unknown parameters).
In what follows we suppose that the bounds $L$
and $R_0$ are known {\em a priori}.
\par
We analyze here the following adaptive version of Algorithm \ralgo{2} (
we leave the construction and analysis
of adaptive version of Algorithm \ralgo{22} as an exercise to the reader):
\balgo{adaptive}
\begin{description}
\item[Initialization:] Set $y_0=x_0$,
$m=\left[{1\over 2}\log_2{\mu(d)N\over
A(d)\log_2N}\right]-1$
\footnote{$^)$ here $[a]$ stands here for the largest integer less or equal to $a$}$^)$,
$N_0=[ N/m]$,
and
$$
R_k=2^{-k}R_0, \quad k=1,... ,m.
$$
\item[Stage] $k=1,...,m$: Compute $y_k = \wh{x}_{N_0}(y_{k-1},R_{k-1})$ with
$\gamma_k = {LR_{k-1} \over \sqrt{2 \mu(d) A(d) }}$.
\item [Output:] $\wh{x}_N=\argmin_{k=1,...,m}
f(y_k)$.
 \end{description}
\ealgo
\bthm  The approximate solution $\wh{x}_N$ satisfies for $N\ge 4$
\bse
f(\wh{x}_N)-f^* & \le & 2\left({16L^2 A(d)\log_2N\over \mu(f)^{2\over\rho }\mu(d)N
 }\right)^{\rho\over 2(\rho-1)}.
\ese
\ethm{adapt}
\section{Generating dual solutions}\label{sc-PD}

In order to speak about primal-dual solutions, we need to
fix somehow the structure of objective function in problem
(\ref{prob-Main}). Let us assume that
$$
f(x) \; = \; \max\limits_{w \in S} \; \Psi(x,w), \quad x
\in Q,
$$
where $S$ is a closed convex set, and function $\Psi$ is
convex in the first argument $x \in Q$ and concave in the
second argument $u \in S$. Let us assume that $\Psi$ is
subdifferentiable in $x$ at any $(x,w) \in Q \times S$.
Then we can take
\beq\label{eq-UX}
\ba{rcl}
f'(x) & = & \Psi'_x(x,w(x)),\\
\\
w(x) & \in & \mbox{Arg}\max\limits_{w\in S} \Psi(x,w).
\ea
\eeq
Thus, we can define the dual function $\eta(w) \; = \;
\min\limits_{x \in Q} \; \Psi(x,w)$, and the dual
maximization problem
\[
\mbox{Find $f^*$} =  \max\limits_w \{ \eta(w):\; w \in
S \}.
\]
For any $w \in S$, we assume that $\Psi(\cdot,w)$ is
uniformly convex on $Q$ with convexity parameters $\rho=\rho(\Psi)$ and $\mu(\Psi))$.
\par
Let for $x,w\in \mathbb R^n$ and let \[
\Psi(x,w)=\langle w,x\rangle-\half\|w\|_q^2, \;\;\;2\le q<\infty.
\]
Clearly, $\Psi(x,w)$ is convex in $x$ and concave in $w$. Further,
\[f(x)=\max_{x\in \mathbb R^n}\Psi(x,w)=\half \|x\|_p^2, \;\;p={q\over q-1},
 \]is strongly convex with respect to $\|\cdot\|_1$ on $\mathbb R^n$ with $\mu(f)=O(1)n^{p-1\over p}$, $f'(x)=w(x)$, where
\[
w^{(i)}(x)=\|x\|_p^{q-2\over q-1}|x^{(i)}|^{1\over q-1} \sign(x^{(i)}).
\]
\bthm
Let assumptions of Theorem \rth{th-Mult1} hold and let $\wh{x}_\e(y_0,R_0)$ be the approximate solution,
supplied by Algorithm \ralgo{1}.
Define $\bar w_{N_m} = {1 \over 1+N_m}
\sum\limits_{i=0}^{N_m} w(x_i)$. Then
\[
f(\wh{x}_\e(y_0,R_0)) - \eta(\bar
w_{N_m}) \leq C(\rho)\, \epsilon,
\]
where
\[
C(\rho)\leq1+3\,{6^{1\over \rho-1}+2^{1\over \rho}\rho^{1\over \rho-1}\over \rho^{\rho\over \rho-1}}
+
{6\over 2^{\rho-1\over \rho}\rho}.
\]
Furthermore, when the objective $f$ is strongly convex ($\rho=2$),
\[
f(\wh{x}_\e(y_0,R_0)) - \eta(\bar
w_{N_m}) \leq 8.5\,\e.
\]
\ethm{th-PD}
%
%

\section{Stochastic programming with uniformly convex objective}
\label{stoch}
In order to rewrite the results of sections \ref{deterministic} in the stochastic
framework we substitute for $f'(x_k)$ its observation $g_k=f'(x_k)+\xi_k$
into the iteration of Algorithm \ralgo{DA-met}.
The following statement is a stochastic counterpart of
Proposition \rp{DA}:
\bpro
 Let $x_k$, $k=0,1,...$ be the search points of Algorithm \ralgo{DA-met} with $g_k$ substituted
 for $f'(x_k)$. Then for any $x\in Q\cap B_R(\bar{x})$,
\begin{equation}\label{eq:eq-DAstoch1}
\sum_{i=1}^k \lambda_i\la f'(x_i), x_i-x\ra\leq
d_{\bar{x},R}(x)
\beta_{k+1} +{R^2
\over 2 \mu(d)} \sum\limits_{i=0}^k {\lambda_i^2 \over
\beta_i} \| f'(x_i)\|^2_*+
\sum\limits_{i=0}^k \zeta_i,
\end{equation}
where
\be
\|\zeta_i\|_*\le 2\lambda_i\|\xi_i\|_*R,\;\;\;\zeta_i\le-\lambda_i\la \xi_i, \tilde{x}_i-x \ra
+{R^2\lambda_i^2 \|\xi_i\|_*^2
\over 2 \mu(d)\beta_i},
\ee{eq-DAstoch2}
and $(\tilde{x}_i), i=1,...,k$ are $\cF_{i-1}$-measurable
random vectors, $\tilde{x}_i\in Q\cap B_R(\bar{x})$
\epro{stochP}
In this section we propose two families of multi-stage methods for uniformly
convex stochastic programming problem described in section \ref{prstat}.
The first one is based on the dual averaging scheme with the prox-function which satisfies the condition
\rf{quadgr} of quadratic growth. As we have already mentioned, one can easily obtain the bounds
for the average value of the objective at the approximate solution, generated by the stochastic counterpart of Algorithm \ralgo{5} and \ralgo{22}. On the other hand, the methods derived from those,
presented in section \ref{sec:smallballs}, better suit the case when the confidence bounds on
the error of the approximate solutions are required.
\par
\subsection{Expectation bounds for methods with prox-function of quadratic growth}
\label{sec:expect}When taking the expectation with respect to the distribution of $\xi_i$ we obtain the following
simple counterpart of Lemma \ref{lm-DA}:
\BL\label{lem-estoch}
Let  $\bar x \in Q$ satisfy
$\bE\| \bar x - x^* \|^2\le  R^2$, where $x^*$ is the optimal
solution of problem (\ref{prob-Aux}), and let $ \lambda_k=1$ and $\beta_k=\g\sqrt{N+1}$, $k=0,..,N$.
Suppose that
Assumptions \ras{stoch1} and \ras{stoch2} hold.
Then the
approximate solution supplied by Algorithm \ralgo{DA-met}
satisfies
\bse
\bE f(x_N(\bar{x},R))-f^*&\le& {1\over N+1}\sum_{i=0}^N  \bE\la
f'(x_i), x_i-x^*\ra \nn
&\leq& {1\over \sqrt{N+1}} \left(\g \bE d_{\bar{x},R}(x^*)
+{R^2 (L^2+\sigma^2)\over 2\mu(d)\gamma}\right),\\
\bE\|x_N(\bar{x},R)-x^*\|^\rho  &\leq &\; {1\over
\mu(f)\sqrt{N+1}} \left(\g \bE d_{\bar{x},R}(x^*) +{R^2
(L^2+\sigma^2) \over 2\mu(d)\gamma}\right).
\ese
\EL
Suppose now that $\bE\|\bar{x}-x^*\|^2\le r^2$. Using the relation
$
d_{\bar{x},R}(x^*)\le C(d){r^2\over R^2}
$
we get the following (cf Corollary \rc{cor-Stage2})
\BC \label{the:exp-del}
Suppose
that $\bar{x}\in Q$ satisfy
\[
\bE\|\bar{x}-x^*\|^2\le r^2,
\]
and let
\[
\g={R^2\over r}\sqrt{L^2+\sigma^2\over 2C(d)\mu(d)},
\]
Then
\be
\label{eq:cdf}
\bE f(x_N(\bar{x},R))-f^*&\le& r\sqrt{2C(d)(L^2+\sigma^2)\over \mu(d)(N+1)},\\
\bE\|x_N(\bar{x},R)-x^*\|^\rho&\le& {r\over
\mu(f)}\sqrt{2C(d)(L^2+\sigma^2)\over \mu(d)(N+1)}.
\ee{cdr}
\EC
%
%
When comparing the above statement to the result of Corollary \rc{cor-Stage2}
 we observe that the only difference between the two is that
in Corollary \rc{exp-del} the quantity $L^2$ is substituted with
$L^2+\sigma^2$. When modifying in the same way the parameters of Algorithm \ralgo{22}
we obtain the multistage procedure for the stochastic problem.
\par
Assume
that the parameters $L$, $\rho,\,\mu(f)$ and the diameter $R_0$ of
$Q$ are known.  The method starts from an arbitrary $x_0 \in Q$.
\balgo{5s}
\begin{description}
\item[Initialization:]
Set $y_0 = x_0$, $\tau={2(\rho-1)\over \rho}$ and $m =\lfloor\log_2 {\mu(f) \over
\epsilon} R_0^\rho\rfloor+1$.
\item[Stage] $k = 1, \dots , m$:
\begin{enumerate}
\item Define $N_k = \lfloor 2^{\tau k} {4 (L^2+\sigma^2) C(d)
\over \mu^2(f)
\mu(d) R_0^{2(\rho-1)}}\rfloor$ and $r^\rho_k = 2^{-k}R^\rho_0$.
\item Compute $y_k = x_{N_k}(y_{k-1},R_0)$ with
$\gamma_k = {R_0^2\over r_{k-1}}\sqrt{L^2+\sigma^2\over 2C(d)\mu(d)}$.
\end{enumerate}
\item[Output:] Set the approximate solution $\wh{x}_\e=y_m$.
\end{description}
\ealgo
 We have
the following stochastic analogue of Theorem
\rth{multi-qdet}:
\bthm Suppose that
\[
N\ge N(\e)={4^{\tau+1} (L^2+\sigma^2) C(d)\over \mu(f)^{2\over \rho}\mu(d)}
\e^{-\tau}.
\]
Then the approximate solution $\wh{x}_N$, provided by
Algorithm \ralgo{5s} satisfies:
\[
\bE f(\wh{x}_\e)-f^*\le \e.
\]
\ethm{multi-qstoch}
 The proof of the theorem follows the
lines of that of Theorem \rth{multi-qdet}. It suffices
to substitute the bounds \rf{cdf} and \rf{cdr} for those
of \rf{eq-FB0} and \rf{eq-NB0}. We leave this simple
exercise to the reader.
\par
The method can be rewritten for the case when
the total number $N$ of calls to the oracle is fixed.
\par
Suppose that
\[N\ge 2^\tau(2^\tau+1) {4(L^2+\sigma^2)C(d)\over
\mu^2(f)\mu(d) R_0^{2(\rho-1)}}.
\] Consider the following procedure:
\balgo{22s}
\begin{description}
\item[Initialization:] Set $y_0 = x_0$, $\tau={2(\rho-1)\over \rho}$, compute
$N_j=\left\lfloor 2^{\tau j}{4(L^2+\sigma^2)C(d)\over
\mu(f)^2\mu(d)R_0^{2(\rho-1)}}\right\rfloor$, while $\sum_{j} N_j\le
N$. Set $m(N)=\max\{k:\;\sum_{j=1}^k N_j\le N\}$.
\item[Stage] $k = 1, \dots , m(N)$:\\
Set $r_k^\rho=2^{-k}R^\rho_0$. Compute $y_k =
x_{N_k}(y_{k-1},R_0)$ with
$\gamma_k = {R_0^2\over r_{k-1}}\sqrt{L^2+\sigma^2\over 2C(d)\mu(d)}$.

\item[Termination:] Set the approximate solution $\wh{x}_N=y_{m(N)}$.
\end{description}
\ealgo
\BC\label{the:mult2-sto} We have
\bse
 \bE f(\wh{x}_N) - f^* \leq 2\left({8(L^2+\sigma^2)C(d)\over \mu(f)^{2\over \rho}\mu(d)N}\right)^{1/\tau}.
\ese
\EC
\par
Exactly in the same way it was done in the deterministic settings, we can
provide an adaptive version of the method. To
this end the adaptive method of Algorithm \ralgo{adaptive} for deterministic problem
should be slightly
modified: we have to change the way the
approximate solution $\wh{x}_N$ is formed, as the exact
observations of the objective function are not available
anymore. Fortunately, we can take as the output of the algorithm the approximate solution $y_m$,
generated at the last stage.
\par
Consider the following procedure:
\balgo{stoca}
\begin{description}
\item[Initialization:] Set $y_0=x_0$,
$m=\left[{1\over 2}\log_2{\mu(d)N\over
C(d)\log_2N}\right]-1$, $N_0=[N/m]$,
$r_k=2^{-k}R_0$, $k=1,... ,m$.
\item[Stage] $k=1,...,m$: Compute $y_k = x_{N_0}(y_{k-1},R_0)$ with
with
%
$\gamma_k = {R_0^2\over r_{k-1}}\sqrt{L^2+\sigma^2\over 2C(d)\mu(d)}$.

\item [Termination:] Set the approximate solution $\wh{x}_N=y_m$.
 \end{description}
\ealgo
\bthm The approximate solution $\wh{x}_N$, supplied by
Algorithm \ralgo{stoca}, satisfies for $N>4$:
\[
\bE f(\wh{x}_N)-f^* \le 4\left({16(L^2+\sigma^2) C(d)\log_2N\over \mu(f)^{2\over\rho }\mu(d)N
 }\right)^{\rho\over 2(\rho-1)}.
\]
\ethm{stocha}

\subsection{Confidence sets for uniformly convex stochastic programs}
In this section we
establish confidence bounds for the approximate solutions, delivered by multistage stochastic algorithms. Consider dual averaging
Algorithm \ralgo{DA-met} in which we substitute the exact subgradient with the observation $g_k=f'(x_k)+\xi_k$.
Let $\delta_N(\bar{x},R)$ be the gap value, defined in \rf{gapv}.
\bpro
Let $\bar{x}$ be a  point of $Q$,
$ \lambda_k=1$ and $\beta_k=\g\sqrt{N+1}$, $k=0,..,N$.
Suppose that Assumptions \ras{stoch1}--\ras{stoch3} hold.
Then
\begin{equation}\label{eq:simplmdp}
\Prob_{\bar{x}}\left[\delta_N(\bar{x},R)\ge {1\over \sqrt{N+1}}\left(\g A(d)+{R^2(L^2+\sigma^2)\over 2\g\mu(d)}\right)+ 2R\sigma \sqrt{3\ln\al^{-1}\over N+1}\right]\le \alpha.
\end{equation}
\epro{mdp}
From \rf{simplmdp} we obtain immediately:
\BC\label{the:Stage1-mdp}
Let $\bar{x}$ be a  point of $Q$. Let
\[
\g  =  R \sqrt{L^2+\sigma^2\over  2\mu(d) A(d)}.
\]
Then for all $\al\ge 0$, the approximate solution $x_N(\bar x, R)$ of Algorithm \ralgo{DA-met}
satisfies
\begin{equation}\label{eq:NB1-mdp}
\Prob_{\bar{x}}\left[\delta_N(\bar{x},R)\le 2R \left[\sqrt{ A(d)(L^2+\sigma^2) \over 2\mu(d) (N+1)}+
\sigma \sqrt{\ln 3\al^{-1}\over N+1}\right]\right]\ge 1-\alpha.
\end{equation}
\EC
%
\par
Corollary \rc{Stage1-mdp} allows us to compute the confidence sets for approximate solutions, provided by stochastic analogues of Algorithms \ralgo{1} and \ralgo{2} exactly in the same way as it was done in section \ref{sec:smallballs}. For the sake of conciseness we present here only the result for the setting
when the total number $N$
of subgradient observations  is fixed and the convexity parameters of the objective are unknown.
\balgo{adaptive-s}
\begin{description}
\item[Initialization:] Set $y_0=x_0$,
$m=\left[{1\over 2}\log_2{\mu(d)N\over
A(d)\log_2N}\right]-1$, $N_0=[N/m]$, and
$$
R_k=2^{-k}R_0, \quad k=1,... ,m.
$$
\item[Stage] $k=1,...,m$: Compute $y_k = \wh{x}_{N_0}(y_{k-1},R_{k-1})$ with
$\gamma_k = {R_{k-1} \sqrt{N_0(L^2+\sigma^2) \over 2\mu(d) A(d) }}$.
\item [Output:] $\wh{x}_N=y_m$.
 \end{description}
\ealgo
\bthm  Let $\al\ge 0$.
Then
the approximate solution $\wh{x}_N$ satisfies for $N\ge 4$
\bse
\Prob\left[f(\wh{x}_N)-f^* \le \e(N,\alpha)\right]\ge 1-\alpha,
\ese
where
\bse
\e(N,\alpha)=4\left({16\over (N_0+1)\mu(f)^{2\over \rho}}\right)^{\rho\over 2(\rho-1)}
\left(
\sqrt{(L^2+\sigma^2) A(d) \over 2\mu(d)}
+
\sigma \sqrt{3\ln \left({\log_2 N\over 2\al}\right)}
\right)^{\rho\over \rho-1}.
\ese
\ethm{adapt-s}
\section{Computational issues}\label{sec-comp}
The interest of the proposed algorithmic schemes is conditioned
by our ability to compute efficiently the optimal solution $ \pi_{z, R, \beta}(s)$
of the optimization problem
\rf{V}.
We present here two important examples in which the problem \rf{V}
can be solved quite efficiently. These are the standard simplex and the hyperoctahedron settings.
\par
Let us measure the
distances in $E = \mathbb R^n$ in $l_1$-norm:
$$
\ba{c}
\| x \| = \| x \|_1 = \sum\limits_{i=1}^n | x^{(i)} |.
\ea
$$

\subsection{Simplex setup}
Let $n\ge 2$ and let
\[
Q=\{x\in \mathbb R^n|\,x\ge 0, \|x\|_1=1\}
\]
be the standard simplex. We are to show how the problem
\rf{V} can be solved in this case. The problem \rf{V} on
$Q_{R}(z)$ for the function $d$ as in \rf{DistL1} writes
\bse\begin{array}{l}
\min_{x,u,v} \left\{\sum_{i=1}^n  \,[s_ix_i+u_i\ln u_i+v_i\ln v_i]:\;\;\;\sum_{i=1}^n [u_i+v_i]=R,\;\;\;\sum_{i=1}^n x_i=1,\right.\\
\left.x_i=z_i+u_i-v_i, \;\;u_i\ge 0, \;\;v_i\ge 0,\;\;\;x_i\ge 0,\;\;i=1,...,n.\right\}\end{array}
\ese
When eliminating the ``$x$'' variable and dualizing the coupling constraints we obtain the equivalent problem
\be
&&\max_{\lambda,\mu}\left\{\underline{L}(\lambda,\mu)\equiv \min_{u,v} L(u,v,\lambda,\mu)\,:
\;\;z_i+u_i-v_i\ge 0,\;\;i=1,...,n\right\},\\
&&\mbox{where}\nn
&&\begin{array}{l}
L(u,v,\lambda,\mu)=\sum_{i=1}^n[r_iv_i+t_iu_i+u_i\ln u_i+v_i\ln v_i]-\lambda R-\mu\,:\\
r_i=s_i+\lambda-\mu
,\;\;t_i=-s_i+\lambda+\mu.\nn
\end{array}
\ee{V-dist}
The dual problem \rf{V-dist} can be solved using a conventional method of convex optimization (ellipsoid or level),
given the solution of the problem
\[
\min_{u,v} \left\{L(u,v,\lambda,\mu)\,:
\;\;z_i+u_i-v_i\ge 0,\;\;i=1,...,n\right\}.
\]
Note that the latter problem an be decomposed into $n$ 2-dimensional problems
\be
\min_{u,v}\; su+tv+u\ln u+v\ln v, \;\;u\ge v-z.
\ee{62}
One way to compute the minimizer is to compute  the solution $(\bar{u},\bar{v})$ to the problem
\bse
\min_{u,v}\; [\psi(u,v)=su+tv+u\ln u+v\ln v], \;\;u=v-z,
\ese
namely,
\[
\bar{u}={1\over 2}\left(\sqrt{z^2+4e^{-2-s-t}}-z\right), \;\;\bar{v}={1\over 2}\left(\sqrt{z^2+4e^{-2-s-t}}+z\right)
\]
and to see if the subgradient
\[
\psi'(u,v)=\left(\begin{array}{c} s+\ln u+1\\
t+\ln v+1\end{array}\right).
\]
satisfies
\[
\psi'_u(\bar{u},\bar{v})+\psi'_v(\bar{u},\bar{v})=0\;\;\; \mbox{and}
\;\;\;\psi'_u(\bar{u},\bar{v})-\psi'_v(\bar{u},\bar{v})>0.
\]
If this is the case, we take $\bar{u},\bar{v}$ as the minimizers, if not, the inequality constraint is not active at the optimal solution of \rf{62} and we take
\[
\bar{u}=e^{-1-s},\;\;\bar{v}=e^{-1-t}.
\]
\subsection{Hyperoctahedron setup}
Let now $Q$ be a standard hyperoctahedron: $Q=\{x\in
\mathbb R^n|\,\|x\|_1\le 1\}$. Let us see how the solution
to \rf{V} can be computed in this case.
\par
When writing
\[
x_i=w_i-y_i, \;\;w_i,\;y_i\ge 0, \;\;\sum_{i=1}^n [w_i+y_i]=1,
\]
the problem \rf{V} on $Q_{R}(z)$ can be rewritten as
\bse\begin{array}{l}
\min\limits_{w,y,u,v} \left\{\sum_{i=1}^n
\,[s_i(w_i-y_i)+u_i\ln u_i+v_i\ln v_i]:\;\;\;
\sum_{i=1}^n [u_i+v_i]=R,\;\;\;\sum_{i=1}^n [w_i+y_i]=1,\right.\\
\left.w_i-y_i=z_i+u_i-v_i, \;\;u_i\ge 0, \;\;v_i\ge 0,\;\;w_i\ge 0,\;\;y_i\ge 0,\;\;i=1,...,n.\right\}\end{array}
\ese
When dualizing the coupling constraints we come to
\bse
&&\max_{\lambda,\mu}\left\{\underline{L}(\lambda,\mu)\equiv \min_{u,v,w,y} L(u,v,w,y,\lambda,\mu)\,:\right.\\
&&\left.\;\;z_i+u_i-v_i-w_i+y_i= 0,\;\;w_i\ge 0, \;\;y_i\ge 0,\;\;i=1,...,n\right\}
\ese
where
\bse
\begin{array}{l}
L(u,v,w,y,\lambda,\mu)=\sum_{i=1}^n[r_iv_i+t_iu_i +\mu (w_i+y_i)+u_i\ln u_i+v_i\ln v_i]-\lambda R-\mu\,:\\
r_i=s_i+\lambda
,\;\;t_i=-s_i+\lambda.\nn
\end{array}
\ese
The computation of the dual function $\underline{L}(\lambda,\mu)$ boils down to evaluating solutions
to $n$ subproblems
\be\begin{array}{l}
\min_{u,v}\; su+tv+\lambda(w+y)+u\ln u+v\ln v, \\
z+u-v-w+y=0,\;\; w\ge 0,\;\;y\ge 0.
\end{array}
\ee{small2}
It is obvious that either $w$ or $y$ vanishes, and to find the solution to \rf{small2} it suffices to compare
the optimal values of the problems
\bse
\min_{u,v}\; \psi_w(u,v)=su+tv+\lambda (z+u-v)+u\ln u+v\ln v, \;\;z+u-v\ge 0,\;\;(\mbox{case $y=0$}),\\
\min_{u,v}\; \psi_y(u,v)=su+tv-\lambda (z+u-v)+u\ln u+v\ln v, \;\;z+u-v\le 0,\;\;(\mbox{case $w=0$}),
\ese
which are the same problems as \rf{62} in the previous section.
\section*{Acknowledgements} {\em The authors would like to acknowledge insightful and motivating comments of Prof. Peter Glynn, which were extremely helpful to them upon completion of this paper.}

\appendix
\section{\protect{Lower complexity bound for uniformly convex optimization}}\label{app-lower}
 For the sake of simplicity we consider here the minimization problem
 \beq\label{prob-Main-lower}
\min\limits_{x} \{ f(x):\; x \in Q \},
\eeq
 over the domain $Q$ which is an Euclidean ball:
 \[
 Q=\{x\in \mathbb R^n|\,\|x\|_2\le R\}.
  \]
The lower bound below can be reproduced for domains of different geometry when following the construction in \cite[chapter 3]{nemyud:83}.
\par
Let $\cF_R(L, \rho)$ be a class of Lipschitz continuous and uniformly convex functions on $Q$, with Lispchitz constant $L$ and parameters of uniform convexity $\rho$ and $\mu(\rho)=1$, when measured with respect to the Euclidean norm. Note that each problem (\ref{prob-Main-lower}) from the class is solvable; we denote $f^*$ corresponding optimal value.
\par
We equip $\cF_R(L, \rho)$ with a first order oracle and define the {\em analytical complexity} ${\cal A}(\epsilon)$ of the class in the standard way:
\[
{\cal A}(\epsilon)=\inf_{\cal M} {\cal A}(\epsilon,{\cal M});
\]
where the (analytical) complexity ${\cal A}(\epsilon,{\cal M})$ of a method $\cal M$ is the minimal number of oracle calls (steps of $\cal M$) required by $\cal M$ to solve any problem of the class $\cF_R(L, \rho)$ to absolute accuracy $\epsilon$ -- find an approximate solution $\bar{x}$ such that $f(\bar{x})-f^*\le \epsilon$.
\bthm
Assume that\footnote{Note that any uniformly convex, with parameters $\rho$ and $\mu=1$, function $f$ on $Q$ clearly satisfies $L\ge C(\rho)R^{\rho-1}$, cf. (\ref{def-Conv12}). }
\beq\label{Lbound}
L\ge 2^{\rho-2}\rho R^{\rho-1}.
\eeq Then the analytical complexity ${\cal A}(\epsilon)$ of the class $\cF_R(L, \rho)$ admits the lower bound:
\[
{\cal A}(\epsilon)\ge \min\left\{n, \left\lfloor {L^2R^2\over 16\epsilon^{2}}\right\rfloor ,\,\left\lfloor L^{2}\over 8\epsilon^{{2(\rho-1)\over \rho}}\right\rfloor\right\}
\]
(here $\lfloor \cdot\rfloor$ stands for the integer part).
\ethm{th-lower}
\pr The proof of the lower bound reproduces the standard
reasoning of \cite[chapter 3]{nemyud:83}.
It suffices to prove that if
$\epsilon\in (0,1)$ is such that
\[
M=\left\lfloor\min\left\{{L^2R^2\over 16\epsilon^{2}},\
{L^{2}\over 8\epsilon^{{2(\rho-1)\over \rho}}}
\right\}-0\right\rfloor\le n
\]
then
the complexity ${\cal A}(\varepsilon)$ is at least $M$. Assume
that this is not the case, so that there exists a method ${\cal
M}$ which solves all problems from the family in question in no
more than $M-1$ steps. We assume that ${\cal M}$ solves any
problem exactly in $M$ steps, and the result always is the last
search point. Let us set
\beq\label{iniplow}
\delta=\min\left\{
{LR\over 4\sqrt{M}},\,
{L^{\rho\over \rho-1}\over 8M^{\rho\over 2(\rho-1)}}
\right\} -
\epsilon,
\eeq so that $\delta > 0$ by definition of $M$.
Now for $\lambda>0$
consider the family ${\cal F}_0$ comprised of functions
\[ f(x) =\half L\max_{1\le i\le M} (\xi_ix^i + d_i) + 2^{\rho-3} \|x\|_2^\rho,
\]
where $\xi_i\in \{\pm 1\}$ and
$0<d_i<\delta$, $i=1,...,M$. Note that all functions of the family are well-defined, since
$M\le n$. Furthermore, by (\ref{Lbound}) $f$ is Lipschitz-continuous with Lipschitz constant $\leq L$, and by Lemma 4 of \cite{Nest:2008} the function
$2^{\rho-3}\|x\|_2^\rho$ is uniformly convex with corresponding parameters $\rho$ and $\mu=1$, thus $f(x)$ are  uniformly convex with parameters $\rho$ and  $\mu(f)=1$.
\par
Let us consider the following construction.
Let $x_1$ be the first search point generated by ${\cal M}$; this
point is instance-independent. Let $i_1$ be the index of the
largest in absolute value of the coordinates of $x_1$. We set
$\xi_{i_1}^*$ to be the sign of the coordinate and put $d_{i_1}^* =
\delta/2$. Now let ${\cal F}_1$ be comprised of all functions from
${\cal F}$ with $\xi_{i_1}= \xi_{i_1}^*$, $d_{i_1}=d_{i_1}^*$ and
$d_i\le\delta/4$ for all $i\ne i_1$. It is clear that all the
functions of the family ${\cal F}_1$ possess the same local
behavior at $x_1$ and are positive at this point.
\par
Now let at the step $k+1$ $i_{k+1}$ be the index of largest in absolute value of the coordinates of
$x_{k+1}$ with indices different from $i_1,...,i_k$. We
define $\xi_{i_{k+1}}^*$ as the sign of the coordinate, put $d_{i_{k+1}}^*=2^{-(k+1)}\delta$,
and
define ${\cal F}_{k+1}$ as the set of those
functions
from ${\cal F}_k$ for which $\xi_{i_{k+1}} = \xi_{i_{k+1}}^*$,
$d_{i_{k+1}} = d_{i_{k+1}}^*$ and $d_i\le 2^{-(k+2)}$ for
$i$ different from $i_1,...,i_{k+1}$.
\par
It is immediately seen that
the family ${\cal F}_{k+1}$ satisfies the predicate:
\paragraph{
${\cal P}_k$:} {\sl the first $k+1$ points $x_1,...,x_{k+1}$ of the
trajectory of ${\cal M}$ as applied to a function from the
family do not depend on the function, and all the functions from
the family coincide with each other in a certain neighborhood of the
$k+1$-point set $\{x_1,...,x_{k+1}\}$ and are positive at this set.}
\par
Observe that after $M$ steps
we end up with the family ${\cal F}_M$ which consists of exactly one
function
\[
f(x)=\half L\max_{1\le i\le M} (\xi_i^* x_i + d_i^*)+2^{\rho-3}\|x\|_2^\rho
\]
such that $f$ is positive along the sequence
$x_1,...,x_M$ of search points
generated by ${\cal M}$ as applied to the function.
Let now
\[
\overline{x}=-\lambda \sum_{i=1}^M {\xi_i^*}e_i,
\]
where $e_i$ stand for basic orths of ${\Bbb R}^n$, and
\[
\lambda=\min\left\{{R\over \sqrt{M}},\,
\left({2^{2-\rho}L\over \rho M^{\rho/2}}\right)^{1\over \rho-1} \right\},
\]
so that $\overline{x}$ belongs to $Q$. Consider the case of $\lambda=
\left({2^{2-\rho}L\over \rho M^{\rho/2}}\right)^{1\over \rho-1}
$.
We have
\bse
f^* &\le& f(\overline{x}) < -\half L\lambda+
2^{\rho-3}\|\overline{x}\|_2^\rho+\delta=-\half L\lambda+2^{\rho-3}M^{\rho\over 2}\lambda^{\rho}+\delta\\
&\le & -{L^{\rho\over \rho-1} \over M^{\rho\over 2(\rho-1)}}
\left[{2^{2-\rho\over \rho-1}\over 2\rho^{1\over \rho-1}}-
{2^{\rho-3}2^{\rho(2-\rho)\over \rho-1}\over \rho^{\rho\over\rho-1}}\right]+\delta\\
&=& -{L^{\rho\over \rho-1} \over M^{\rho\over 2(\rho-1)}} {2^{2-\rho\over \rho-1}\over 2\rho^{1\over \rho-1}}[1-\rho^{-1}]+\delta\le -{L^{\rho\over \rho-1} \over 8M^{\rho\over 2(\rho-1)}}+\delta\le -\epsilon
\ese
(the concluding inequality follows from (\ref{iniplow})). In the case of $\lambda=R/\sqrt{M}$ we have
\[
f^*\le f(\overline{x})\le -\half L\lambda+
2^{\rho-3}\|\overline{x}\|_2^\rho+\delta\le -{LR\over 4\sqrt{M}}+\delta\le -\epsilon.
\]
Thus, in both cases we have
$
f(x_M)-f^* > 0-(-\epsilon)=\epsilon.
$
Since, by construction, $x_M$ is the result obtained by ${\cal M}$
as applied to $f$, we conclude that ${\cal M}$ does not solve the
problem $f$ within relative accuracy $\varepsilon$, which is the
desired contradiction with the origin of $M$.\qed

\section{Proofs}\label{app-proofs}
\subsection{Proof of Lemma \ref{lm-Symm}}
Consider two points $x_i \in Q^0$, $i = 1,2$. Suppose that
$$
\ba{c}
x_i = u_i - v_i, \quad u_i \in \alpha_i Q, \quad v_i \in
(1 -
\alpha_i)Q,\quad \alpha_i \in [0,1],\\
\\
f^0(x_i) = f(u_i) + f(v_i), \quad i = 1, 2.
\ea
$$
Let us choose an arbitrary $\alpha \in [0,1]$. Then,
$$
\ba{rcl}
x(\beta) & \Def & \beta x_1 + (1 - \beta) x_2 \\
\\
& = & \beta (u_1 - v_1) + (1 - \beta) (u_2 - v_2) \\
\\
& = & \beta u_1 + (1 - \beta) u_2 - (\beta v_1 +
(1-\beta)v_2).
\ea
$$
Denote $\gamma = \beta \alpha_1 + (1-\beta)\alpha_2$. Then
$$
1 - \gamma = \beta (1 - \alpha_1) + (1- \beta)(1 -
\alpha_2).
$$
Note that $u_i = \alpha_i \bar u_i$, and $v_i = (1 -
\alpha_i) \bar v_i$ for some $\bar u_i$ and $\bar v_i$
from $Q$, $i = 1, 2$. Therefore, denoting
$$
\tau = \beta \alpha_1/ \gamma, \quad \xi = \beta (1 -
\alpha_1) / (1 - \gamma),
$$
we obtain
$$
\ba{rcl}
x(\beta) & = & \beta \alpha_1 \bar u_1 + (1 - \beta)
\alpha_2 \bar u_2 - (\beta (1 - \alpha_1) \bar v_1 +
(1-\beta)(1 - \alpha_2) \bar v_2)\\
\\
& = & \gamma (\tau \bar u_1 + (1 - \tau) \bar u_2) - (1 -
\gamma)(\xi \bar v_1 + (1 - \xi) \bar v_2) \\
\\
& \Def & \gamma \bar u_3 - (1 - \gamma) \bar v_3
\ea
$$
with some $\bar u_3$ and $\bar v_3$ from $Q$. Hence, $u_3
= \gamma \bar u_3 \in \gamma Q$, and $v_3 = (1 - \gamma)
\bar v_3 \in (1- \gamma) Q$. Consequently, by definition
of function $f^0$ and using inclusions $u_i, v_i \in Q$,
$i=1,2$, we obtain
$$
\ba{rcl}
f^0(x(\beta)) & \leq & f(u_3) + f(v_3) \\
\\
& = & f(\beta u_1 + (1 - \beta) u_2) + f( \beta v_1 +
(1-\beta)v_2)\\
\\
& \leq & \beta f(u_1) + (1 - \beta) f(u_2) - \half \mu
\beta(1 - \beta) \| u_1 - u_2 \|^2 \\
\\
& & + \beta f(v_1) + (1 - \beta) f(v_2) - \half \mu
\beta(1 - \beta) \| v_1 - v_2 \|^2\\
\\
& = & \beta f^0(x_1) + (1 - \beta) f^0(x_2) - \half \mu
\beta(1 - \beta) \left[\; \| u_1 - u_2 \|^2 + \| v_1 - v_2
\|^2 \; \right].
\ea
$$
It remains to note that
$$
\ba{rcl}
2 \| u_1 - u_2 \|^2 + 2 \| v_1 - v_2 \|^2 & \geq & \| u_1
- u_2 - (v_1 - v_2) \|^2 \; = \; \| x_1 - x_2 \|^2.
\ea
$$
\qed
\subsection{Proof of Lemma \ref{lm-DA}}
In view of conditions of the lemma, $x^*
\in Q_{R}(\bar x)$. From the assumptions
on function $f$, we conclude that
$$
\ba{rcl}
\la f'(x_i), x_i - x^* \ra & \geq & f(x_i) - f(x^*),\\
\\
\la f'(x_i), x_i - x^* \ra &
\stackrel{(\ref{def-Conv12})}{\geq} & \mu(f) \cdot \| x_i
- x^* \|^\rho,\quad i = 0, \dots, N.
\ea
$$
Hence,
\bse
(N+1)\dl_N(\bar{x},R)
&\geq&  \sum\limits_{i=0}^N [f(x_i) -
f(x^*)]  \; \geq \; (N+1) [ f(x_N(\bar x, R)) - f(x^*)],
\\
(N+1)\dl_N(\bar{x},R)&\geq& \; \mu(f)
\sum\limits_{i=0}^N \| x_i
- x^* \|^\rho
\geq \mu(f) (N+1) \| x_N(\bar x, R) - x^* \|^\rho.
\ese
It remains to note that $d_{\bar{x},R}(x)\le A(d)$ for any $x\in Q_R(\bar{x})$ use
the inequality \rf{eq-DADelta}.
\qed
\subsection{Proof of Theorem \rth{th-Mult1}}
Indeed, for $k = 0$, (\ref{eq-Rad}) is valid. Assume it is
valid for some $k \geq 0$. Note that
\[
\ba{rcl}
\sqrt{N_{k+1}+1} & \stackrel{\rf{eq-AuxD}}{\geq} &
\left( \left({2^{ k}\over R^{\rho}_0}\right)^{\tau}{8 L^2 A(d)\over \mu^2(f) \mu(d)
}\right)^{1/2} \; = \; 2{ L \sqrt{2 A(d)}\over \mu(f)
R_k^{\rho-1} \sqrt{ \mu(d)}}.
\ea
\]
Therefore, in view of Proposition \rp{DA} and Corollary \rc{cor-Stage1}, we have
$$
\ba{rcl}
\dl_{N_{k+1}}(y_{k},R_{k})& \leq & {L R_k \sqrt{2 A(d)} \over\sqrt{\mu(d)(N_{k+1}+1)}} \le
{\mu(f)\over 2}R^{\rho}_k=\mu(f)R^\rho_{k+1},
\ea
$$
and this is (\ref{eq-Val}) for the next value of the
iteration counter. Further,
\[
\| y_{k+1} - x^* \|^\rho\le
\mu(f)^{-1}\dl_{N_{k+1}}(y_{k},R_{k})\leq \; R_{k+1}^\rho,
\]
and this is (\ref{eq-Rad}) for  $k+1$.

Finally, at the end of the $m$-th stage, in view of
Lemma \ref{lm-DA} and (\ref{eq-Val}) we have
$$
\ba{rcl}
f(\wh{x}_{\epsilon}(y_0,R_0)) - f^* & \leq &
\dl_{N_{m}}(y_{m-1},R_{m-1}) \;
\stackrel{(\ref{eq-Val})}{\leq}\;
 \mu(f) R_m^\rho \; = \; 2^{-m}\mu(f) R_0^\rho \;
 \stackrel{\rf{eq-AuxD}}{\leq} \; \epsilon.
\ea
$$
The complexity of the method can be estimated as follows:
$$
\ba{rcl}
N(\e) & \stackrel{\rf{eq-AuxD}}{\leq} &
\sum\limits_{k=1}^m 2^{k\tau} {4 L^2 A(d)\over \mu^2(f) \mu(d)R_0^{2(\rho-1)}} \;
< \; \left({2^{{m+1}}\over R_0^\rho}\right)^\tau {4 L^2 A(d)\over (2^\tau-1)\mu^2(f) \mu(d)}.
\ea
$$
To conclude (\ref{eq-IT}) it suffices to notice that by \rf{eq-AuxD},
$2^{m+1}\le 4{\mu(f)R_0^\rho\over \e}$.
\qed
\subsection{Proof of Corollary \ref{mult2}}
In the case $N\le \bar{N}$ the corollary follows from the bound of  Corollary
\rc{cor-Stage1} for one-stage method.
When $N\ge \bar{N}$, when following the steps of the proof of Theorem \rth{th-Mult1} we conclude that
\[
f(\wh{x}_N)-f^*\le 2^{-m(N)}\mu(f) R_0^\rho.
\]
Now it suffices to notice that the number $m(N)$ of the stages of the algorithm can be easily bounded: 
\[
\ba{c}
{N\over 2}\; \le \; \sum\limits_{k=1}^{m(N)} N_k \; \le \;
\left({2^{m(N)+1}\over R_0^\rho}\right)^{\tau} {4L^2A(d)\over (2^\tau-1)\mu^2(f)\mu(d)}.
\ea
\]
Thus,
\[
\ba{rcl}
2^{-m(N)} & \le & 2\left({8 L^2A(d)\over \mu^2(f)\mu(d)N}\right)^{1/\tau} R_0^{-\rho},
\ea
\]
and the bound (\ref{eq-it2}) follows.
\qed
\subsection{Proof of Theorem \rth{multi-qdet}}
As in the proof of Theorem
\rth{th-Mult1}, the result of the
theorem follows immediately from the relations:
\be
\| y_k - x^* \|^\rho \leq r^\rho_k = 2^{-k}R^\rho_0
\ee{rad-sto}
and
\be
f(y_k)-f^*\le \mu(f)r^\rho_k\le \mu(f) 2^{-k} R_\rho^2.
\ee{del-sto}
Indeed, using the relations above we write:
\bse
f(\wh{x}) - f^* & \leq &  \mu(f) r_m^\rho=2^{-m}\mu(f)
R_0^\rho\le \e.
\ese
Let us verify the bounds \rf{rad-sto} and \rf{del-sto}.
Assume that \rf{rad-sto} valid for some $k \geq 0$. Note
that
\[
\sqrt{N_{k+1}+1} >  {2^{\tau k/2}\over R_0^{\rho-1}}\left({8 L^2 C(d)\over
\mu^2(f) \mu(d)}\right)^{1/2}= 2{ L\over
\mu(f)r_k 
}\sqrt{2C(d)\over \mu(d)}.
\]
Therefore, in view of  Corollary \rc{exp-del}, we have
\bse
\| y_{k+1} - x^* \|^\rho& \leq & {L r_k \over
\mu(f)}\sqrt{2 C(d) \over \mu(d)(N_{k+1}+1)} \le
{r_k^\rho\over 2}=r_{k+1}^\rho,
\ese
and
\[
f(y_{k+1}) - f^*\le {L r_k}\sqrt{2 C(d) \over
\mu(d)(N_{k+1}+1)}
 \leq {\mu(f)\over 2}r_k^\rho  = \mu(f) r_{k+1}^\rho.
\]
\qed
\subsection{Proof of Theorem \rth{adapt}}
Note that by (\ref{eq-AMu})
 $m$ satisfies $m{\le} {1\over 2}\log_2 {2N
 \over \log_2N}-1\le {1\over 2}\log_2 N$; besides,
\beq\label{eq-Dop1}
\ba{rcl}
2^m & \le & {1\over 2}\sqrt{\mu(d)N\over A(d)\log_2N}.
\ea
\eeq
Assume now that $\mu(f)\le {4L\over R^{\rho-1}_0}\sqrt{A(d)\log_2
N\over \mu(d)N}$. We have
\bse
f(y_1)-f^* & \le & \dl_{N_{0}}(y_{0},R_{0})\le {LR_0}\sqrt{2A(d)\over \mu(d)(N_0+1)}
\; \leq \; {LR_0}\sqrt{2mA(d) \over \mu(d)N}\nn
& \le & {LR_0}\sqrt{A(d)\log_2N\over \mu(d)N}\le
\left({16L^2A(d)\log_2N\over \mu(f)^{2\over \rho}\mu(d)N}\right)^{\rho\over 2(\rho-1)},
\ese
what implies the statement of the theorem in this case. Next, let us
denote $\mu_0=2^{-m}LR_0^{1-\rho}$ so that
\be
 2LR_0^{1-\rho} \sqrt{A(d)\log_2N\over \mu(d)N}\le \mu_0< 4LR_0^{1-\rho}
 \sqrt{A(d)\log_2N\over \mu(d)N},
\ee{mu0bound}
and
$\mu_k=2^{(\rho-1)k}\mu_0$, $k=1,...,m$. Observe that from the
available information we can derive an upper bound on the
unknown parameter $\mu(f)$, namely,
$$
\ba{rcl}
\mu(f) & \le & {L\over R^{\rho-1}_0} \le  \mu_m.
\ea
$$
Suppose now that the true $\mu(f)$ satisfies $\mu_0\le \mu(f)\le \mu_m$. We need
the following auxiliary result.
\BL\label{mlMl}
Let $k^*$ satisfy $\mu_{k^*}\le \mu(f)\le 2^{\rho-1}\mu_{k^*}$. For
$1\le k\le k^*$, the points $\{ y_k \}_{k=1}^m$ generated
by Algorithm \ralgo{adaptive} satisfy the following
relations:
\be\label{eq:rbeforek}
\| y_{k-1} - x^* \| &\leq& R_{k-1} = 2^{-k+1}R_0,\\
\dl_{N_{0}}(y_{k-1},R_{k-1})&\le& \mu_k R_k^\rho= 2^{-k}
\mu_0R_0^\rho.
\ee{dlbeforek}
For $k^*<k\le m$, we have
\be
f(y_k)\le f(y_{k^*})+\mu_{k^*}R_{k^*}^\rho.
\ee{afterk} \EL \proof Let us prove first \rf{rbeforek} and
\rf{dlbeforek}. Indeed, for $k = 1$ \rf{rbeforek} is valid.
 Assume
it is valid for some $k \geq 1$. We write
\bse
\mu(f) & \ge & \mu_{k} \; = \; 2^{(\rho-1)k} \mu_0 \; = \; \left({2^{k}\over R_0}\right)^{\rho-1}
\cdot L  2^{-m}\; \\
&\stackrel{(\ref{eq-Dop1})}{\geq}& \; \left({2^{k}\over R_0}\right)^{\rho-1}\, {2L}\sqrt{A(d)\log_2 N\over \mu(d)N} \; \ge \; {2L\over
R^{\rho-1}_{k}}\sqrt{2A(d)\over \mu(d)(N_0+1)}.
\ese
Therefore,
\beq\label{simplea}
\ba{rcl}
\dl_{N_{0}}(y_{k-1},R_{k-1})&
\stackrel{\rf{eq-FB1}}{\leq} & {L R_{k-1} \sqrt{2 A(d)}
\over\sqrt{\mu(d)(N_{0}+1)}}\; \leq \; \half \mu_k R_k^{\rho-1}
R_{k-1} \; = \; \mu_k R_k^\rho.
\ea
\eeq
That is \rf{dlbeforek}. Moreover,
\[
\ba{rcl}
\| y_{k} - x^* \|^\rho & \le &
\mu(f)^{-1}\dl_{N_{0}}(y_{k-1},R_{k-1})\leq \; {\mu_k\over
\mu(f)}R_{k}^\rho  \le \; R_{k}^\rho,
\ea
\]
and this is \rf{rbeforek} for the next index value. Further, as
in (\ref{simplea}), for $k>k^*$ we have
$$
\ba{rcl}
f(y_{k}) - f(y_{k-1})& \le &\dl_{N_0}(y_{k-1},R_{k-1})\;
\le \; L R_{k-1} \sqrt{2 A(d) \over \mu(d)(N_{0}+1)}\\
\\
& = &2^{k^*-k} LR_{k^*-1}\sqrt{2 A(d)
\over\mu(d)(N_{0}+1)} \; \stackrel{(\ref{simplea})}{\leq}
 2^{k^*-k}\mu_{k^*} R^\rho_{k^*}.
\ea
$$
Then
\[
\ba{rcl}
f(y_{k}) - f(y_{k^*})& =&\sum\limits_{j=k^*+1}^k
f(y_{j}) - f(y_{j-1}) \le \sum\limits_{j=k^*+1}^{k}
2^{k^*-j}\mu_{k^*} R^\rho_{k^*} \le \mu_{k^*}R^\rho_{k^*}.
\ea
\]
This proves the lemma.
\qed
\par
Now we can finish the proof of the theorem. Recall that
$\mu_0\le \mu(f)\le \mu_m$. At the end of the $k^*$-th
stage we have
\bse
f({y}_{k^*}) - f^* & \leq &
\dl_{N_0}(y_{k^*-1},R_{k^*-1})\; {\leq} \;
 \mu_{k^*} R_{k^*}^\rho
 \leq
 2{\mu_{k^*}^{1\over \rho-1}\over \mu(f)^{1\over \rho-1}}\mu_{k^*}R^{\rho}_{k^*}
 \\
 &=&2{\mu_0^{\rho\over \rho-1}R_0^\rho\over \mu(f)^{1\over \rho-1}}\;\stackrel{\rf{mu0bound}}{\le}
 2\left({16L^2 A(d)\log_2N\over \mu(f)^{2\over\rho }\mu(d)N
 }\right)^{\rho\over 2(\rho-1)}.
\ese
\qed
\subsection{Proof of Theorem \rth{th-PD}}
The following result is quite standard (cf. Lemma 3
\cite{BM}).
\BL\label{lm-PD}
Define $\bar x = {1 \over N+1} \sum\limits_{i=0}^N x_i$,
and $\bar w_N = {1 \over N+1} \sum\limits_{i=0}^N w(x_i)$.
Then
\beq\label{eq-PD}
\ba{rcl}
f(\bar x_N) - \eta(\bar w_N) & \leq & l^*_N \Def
\max\limits_x
\left\{ {1 \over N+1} \sum\limits_{i=0}^N \la f'(x_i), x_i
- x \ra - \half \mu(\Psi) \| x - \bar x_N \|^\rho: \; x \in Q
\right\}.
\ea
\eeq
\EL
\proof
Since $\Psi$ is convex in the first argument, for any $x
\in Q$ we have
$$
\ba{rcl}
\la f'(x_i), x_i - x\ra & \stackrel{(\ref{eq-UX})}{=} &
\la
\Psi'_x(x_i,w(x_i)), x_i - x \ra\\
\\
& \geq &\Psi(x_i, w(x_i)) - \Psi(x ,w(x_i)) + \half \mu(\Psi)
\|
x - x_i \|^\rho\\
\\
& = & f(x_i) - \Psi(x, w(x_i)) + \half \mu(\Psi) \| x - x_i
\|^\rho.
\ea
$$
Hence,
$$
\ba{rcl}
l^*_N & = & {1 \over N+1} \max\limits_x \left\{
\sum\limits_{i=0}^N \la f'(x_i),  x_i - x \ra -  \mu(\Psi)
{N+1
\over 2} \| x - \bar x_N \|^\rho: \; x \in Q \right\}\\
\\
& \geq & {1 \over N+1} \max\limits_x
\left\{ \sum\limits_{i=0}^N [\la f'(x_i),  x_i - x \ra -
\half \mu(\Psi) \| x - x_i \|^\rho]: \; x \in Q \right\}\\
\\
& \geq & {1 \over N+1} \max\limits_x
\left\{ \sum\limits_{i=0}^N [ f(x_i) - \Psi(x, w(x_i)) ]: \; x
\in Q \right\} \\
\\
& \geq & f(\bar x_N) - \min\limits_{x \in Q} \Psi(x, \bar
w_N) \; = \; f(\bar x_N) - \eta(\bar w_N).
\ea
$$
\qed

Let us prove now several auxiliary results. Let $l(x)$ be
an affine function on $E$. Let us fix a point $\bar y \in
Q$. Consider the function
$$
\psi(r) = \max\limits_x \{ l(x): \; x \in Q_r(\bar y) \},
\quad r \geq 0.
$$
Note that $\psi(r)$ is an increasing concave function of
$r$ and
$$
\psi(r) \geq \psi(0) = l(\bar y).
$$
Let us fix some $\bar r > 0$ and choose an arbitrary $\bar
x \in Q_{\bar r}(\bar y)$. For some $\mu>0$ define
\beq\label{eq-Lambda}
\lambda^*_{\mu}(x) = \max\limits_y \{ l(y) - \half \mu \|
y - x \|^\rho:\; y \in Q \}.
\eeq
We need to bound from
above the value $\lambda^*_\mu(\bar x)$.
\BL\label{lm-Aux1} For any $b>0$ we have
\be
\lambda^*_\mu(\bar x) \leq \lambda^*_{(1+b)^{1-\rho}\mu}(\bar y)
+  {\mu\over 2b^{\rho-1}}
\bar r^\rho.
\ee{l2}
\EL
\proof
Consider $y_{\mu}(\bar x)$, the optimal solution of
optimization problem in (\ref{eq-Lambda}) with $x =
\bar x$. Then
$$
\lambda^*_{\mu}(\bar x) = l(y_{\mu}(\bar x)) - \half \mu
\| y_{\mu}(\bar x) - \bar x \|^\rho.
$$
On the other hand, for any $b>0$,
$$
\ba{rcl}
\| y_{\mu}(\bar x) - \bar y \|^\rho & \leq & ( \|
y_{\mu}(\bar x) -
\bar x \| + \| \bar x - \bar y \| )^\rho \; \\
&\leq &(1+b)^{\rho-1} \| y_{\mu}(\bar
x) -
\bar x \|^\rho + (1+b^{-1})^{\rho-1} \| \bar x - \bar y \|^\rho\\
\\
& \leq & (1+b)^{\rho-1} \| y_{\mu}(\bar x) - \bar x \|^\rho + (1+b^{-1})^{\rho-1}
\bar r^\rho.
\ea
$$
Hence,
$$
\ba{rcl}
\lambda^*_\mu(\bar x) \leq l(y_{\mu}(\bar x)) - {\mu \over 2}
{\| y_{\mu}(\bar x) -
\bar y \|^\rho\over (1+b)^{\rho-1}} + {1 \over 2b^{\rho-1}} \bar r^\rho \leq
\lambda^*_{(1+b)^{1-\rho}\mu}(\bar y)
+  {\mu\over 2b^{\rho-1}}\bar r^\rho.
\ea
$$
\qed

\BL\label{lm-Aux2}
\[
\lambda^*_{\mu}(\bar y)
\leq \psi(\bar r)+ { \rho-1\over \rho}
\left(2\over \mu\rho\right)^{1\over \rho-1}
\left({\psi(\bar r) -
\psi(0)\over \bar r} \right)^{\rho\over \rho-1}.
\]
\EL
\proof
Indeed, denote $\hat t = \| y_{\mu}(\bar y) -
\bar y \|$. Then
$$
\ba{rcl}
\lambda^*_{\mu}(\bar y) & = & l(y_{\mu}(\bar y)) - \half
\mu \hat t^\rho \; \leq \; \psi(\hat t) - \half \mu
\hat t^\rho \; \leq \; \max\limits_{t \geq 0} \{
\psi(t) - \half \mu t^\rho \} .
\ea
$$
Since $\psi(t)$ is concave,
$$
\psi(t) \; \leq \; \psi(\bar r) + \psi'(\bar r) (t -
\bar r) \; \leq \; \psi(\bar r) + \psi'(\bar r) t.
$$
Note that
\[ \psi'(\bar r)t-\half \mu t^\rho\le  {\rho-1\over \rho}
\left(2\psi'(\bar r)^{\rho}\over \mu\rho\right)^{1\over \rho-1},
\]
thus
\[
\lambda^*_{\mu}(\bar y) \leq \psi(\bar r) +{ \rho-1\over \rho}
\left(2\psi'(\bar r)^{\rho}\over \mu\rho\right)^{1\over \rho-1}.
\] On the other hand,
$$
\psi(0) \leq \psi(\bar r) + \psi'(\bar r)(0 -
\bar r).
$$
Thus, $\psi'(\bar r) \leq {1 \over \bar r} \left(\psi(\bar
r) - \psi(0) \right)$.
\qed
\par
When substituting the result into \rf{l2} we obtain
\BC\label{cor-Aux3}
\beq\label{eq-Aux3}
\ba{rcl}
\lambda^*_{\mu}(\bar x) & \leq & \psi(\bar r) + (1+b){ \rho-1\over \rho}
\left(2\over \mu\rho\right)^{1\over \rho-1}
\left({\psi(\bar r) -
\psi(0)\over \bar r} \right)^{\rho\over \rho-1} +
{\mu\over 2b^{\rho-1}}
\bar r^\rho.
\ea
\eeq
\EC

Let us apply now the above results to Algorithm \ralgo{1}.
Let us choose $\mu = \mu(\Psi)$,
$$
\ba{c}
\bar y = y_{m-1}, \quad \bar x = y_m,\quad \bar r = R_{m-1},
\quad l(x) = {1 \over 1+N_m} \sum\limits_{i=0}^{N_m} \la
f'(x_i), x_i - x \ra,
\ea
$$
where the points $\{ x_i \}_{i=0}^{N_m}$ were generated
during the last $m$th stage of the algorithm. Note that
\beq\label{eq-RInt}
\ba{rcl}
2^m & \geq & {\mu(\Psi) \over \epsilon} R_0^\rho \; \geq \;
2^{m-1}.
\ea
\eeq
Therefore
\beq\label{eq-RR}
\ba{rcl}
{2 \epsilon \over \mu(\Psi)} \; \geq \; \bar r^\rho & = &
2^{1-m}R_0^\rho \; \geq \; {\epsilon \over \mu(\Psi)}.
\ea
\eeq
Further,
$$
\ba{rcl}
\psi(\bar r) & = & \delta_{N_m}(y_{m-1},R_{m-1})\;
\stackrel{(\ref{eq-Val})}{\leq} \; \mu(\Psi) 2^{-m} R_0^\rho \;
\stackrel{(\ref{eq-RInt})}{\leq} \; \epsilon,\\
\\
\psi(0) & = & {1 \over 1 +N_m} \sum\limits_{i=0}^{N_m} \la
f'(x_i), x_i - y_{m-1} \ra \; \geq {1 \over 1 +N_m}
\sum\limits_{i=0}^{N_m} [f(x_i)- f(y_{m-1})]\\
\\
& \geq & f^* - f(y_{m-1}) \;
\stackrel{(\ref{eq-Val})}{\geq} \; - \mu(\Psi)2^{1-m} R_0^\rho
\; \stackrel{(\ref{eq-RR})}{\geq} \; - 2 \epsilon.
\ea
$$
Hence, using the above inequalities in (\ref{eq-Aux3}), we
obtain
\be
\begin{array}{lcl}
\lambda^*_{\mu(\Psi)}(y_m) & \leq &
\e + (1+b){\rho-1\over \rho}
\left(2\over \mu(\Psi)\rho\right)^{1\over \rho-1}
\left({3\e\over \bar r} \right)^{\rho\over \rho-1} +
{\mu(\Psi)\over 2b^{\rho-1}}\bar r^\rho\\
&\leq& \e+{(1+b)(\rho-1)\over 2} \left({6\over \rho}\right)^{\rho\over \rho-1} \kappa^{-{1\over \rho-1}}\e
+{\kappa\e\over 2b^{\rho-1}}\\
&\le &\e\left(1+{(1+b)(\rho-1)\over 2} \left({6\over \rho}\right)^{\rho\over \rho-1}+
{b^{1-\rho}}\right),
\end{array}
\ee{for2}
where we set $\kappa={\mu(\Psi)\bar r^\rho\over \e}$ and used the fact that $1\le \kappa\le 2$ due to (\ref{eq-RR}).

When setting  $b=({\rho\over 6})^{1\over \rho-1}2^{1\over \rho}$ we obtain
\[
\lambda^*_{\mu(\Psi)}(y_m) \leq
\e\left(1+3\,{6^{1\over \rho-1}+2^{1\over \rho}\rho^{1\over \rho-1}\over \rho^{\rho\over \rho-1}}
+
{6\over 2^{\rho-1\over \rho}\rho}\right).
\]
Note that a finer estimate can be obtained for $\rho=2$. To this end it suffices to verify that
for the choice $b=1/3$ the right-hand side of \rf{for2} is decreasing in $\kappa$ for $0\le \kappa\le 2$.
Therefore,
$$
\ba{rcl}
\lambda^*_{\mu(\Psi)}(y_m) & \leq & 8.5 \, \epsilon.
\ea
$$
It remains to note that
$$
\ba{rcl}
\lambda^*_{\mu(\Psi)}(y_m) & = & \max\limits_y \left\{ {1
\over 1 + N_m} \sum\limits_{i=0}^{N_m} \la f'(x_i), x_i -
y \ra - \half \mu(\Psi) \| y - y_m \|^\rho: \; y \in Q \right\},
\ea
$$
and $y_m = x_{N_m}(y_{m-1},R_{m-1})
\stackrel{(\ref{eq-AuxS})}{=} {1 \over 1 + N_m}
\sum\limits_{i=0}^{N_m} x_i$. \qed
\subsection{Proof of Theorem \rth{stocha}}
The proof of the theorem follows
the lines of that of Theorem \rth{adapt}. Using the
notation $k^*$, introduced in Lemma \ref{mlMl}, we get
(cf. \rf{afterk})
\[
\bE f(y_m)\le \bE f(y_{k^*})+\mu_{k^*}r^\rho_{k^*}.
\]
Thus
\[\bE f(y_m)-f^*\le 2\mu_{k^*}r^\rho_{k^*}\le
4\left({16(L^2+\sigma^2) C(d)\log_2N\over \mu(f)^{2\over\rho }\mu(d)N
 }\right)^{\rho\over 2(\rho-1)}.
\]
\qed
\subsection{Proof of Proposition \rp{mdp}}
We need the following result which is essentially known (cf \cite{IbLin}):
\BL\label{boundprobi} Let $\psi_i$, $i=0,...,N$, be Borel functions
on $\Omega$ such that $\psi_i$ is $\cF_i$-measurable,  and let
$\mu_i\geq0$, $\nu_i>0$ be deterministic reals. Assume that for all $i=0,1,2,...$ one has a.s.
\bse
\bE_{i-1}[\psi_i]\leq\mu_i,\;\;\;\bE_{i-1}[\exp\{\psi_i^2/\nu_i^2\}]\leq \exp\{1\},
\ese
 Then for every
$\Lambda\geq 0$
\be
\Prob\left[\sum_{i=0}^N\psi_i>\sum_{i=0}^N\mu_i+\Lambda\sqrt{\sum_{i=0}^N\nu_i^2}\right]
\leq \exp\{-\Lambda^2/3\}\ee{Q.E.Dii}
\EL
For the proof of the lemma see, e.g. section 4.2 of \cite{JNarxiv}.
\par
 Let us return to the proof of the proposition. From \rf{eq-DAstoch2} and Assumption \ras{stoch3} we conclude that
 \[
 \bE_{i-1}\zeta_i\le {R^2\lambda_i^2\sigma^2
\over 2 \mu(d)\beta_i}={R^2\sigma^2
\over 2 \mu(d)\beta_i}
 \](recall that $E_{i-1}\xi_i=0$ and $\tilde{x}_i$ is $\cF_{i-1}$-measurable).
Along with Assumption \ras{stoch3} this implies that random variables $\psi_i=\zeta_i$ satisfy the premises of Lemma \ref{boundprobi}
with $\mu_i={R^2\sigma^2
\over 2 \mu(d)\beta_i}$ and $\nu_i=2R\sigma$.
 Thus by \rf{Q.E.Dii},
\[
\Prob\left[\sum_{i=0}^N\zeta_i\ge {R^2\sigma^2
\over 2 \mu(d)}\sum_{i=0}^N\beta^{-1}_i+2\Lambda R\sigma\sqrt{N+1}\right]\le \exp\{-{\Lambda^2/3}\}\;\;(=\alpha \;\mbox{for $\Lambda=\sqrt{3\ln\al^{-1}}$}).
\]
When substituting $\beta_i=\gamma\sqrt{N+1}$ we conclude \rf{simplmdp} from \rf{eq-DAstoch1}. \qed
%
\subsection{Proof of Theorem \rth{adapt-s}}
Let us denote  $\bar\al={2\al\over\log_2 N}$ and
\bse
a(N_0,\bar \alpha)={2\over \sqrt{N_0+1}}
\left(
\sqrt{(L^2+\sigma^2) A(d) \over 2\mu(d)}
+
\sigma \sqrt{3\ln \bar\al^{-1}}
\right).
\ese
We set
\be
\mu_0=2R_0^{1-\rho}a(N_0,\bar \alpha)\;\;\;
\mbox{and}\;\;\;
\mu_k=2^{(\rho-1)k}\mu_0, \;\;k=1,...,m.
\ee{mui}
Note also that
$$
\ba{rcl}
\mu(f) & \le & {L\over R^{\rho-1}_0},
\ea
$$
and by the definition of $\mu_0$  and $m$ we have $\mu(f)\le\mu_m$.
Suppose first that the true $\mu(f)$ satisfies $\mu_0\le \mu(f)\le \mu_m$. We start with
the following auxiliary result.
\BL\label{mlMl1}
Let $k^*$ satisfy $\mu_{k^*}\le \mu(f)\le 2^{\rho-1}\mu_{k^*}$. Then for
any $1\le k\le k^*$, there exists a set $\cA_k\subset \Omega$ of probability at least $1-k\bar\al$ such that
for $\omega\in \cA_k$ the points $\{ y_k \}_{k=1}^m$ generated
by Algorithm \ralgo{adaptive-s} satisfy
\be\label{eq:rbeforek-s}
\| y_{k-1} - x^* \| &\leq& R_{k-1} = 2^{-k+1}R_0,\\
f(y_k)-f^*&\le& \mu_k R_k^\rho= 2^{-k}
\mu_0R_0^\rho.
\ee{dlbeforek-s}
Further, for $k>k^*$ there is a set $\cC_k\subset\Omega$ of probability at least $1-(k-k^*)\bar\alpha$
such that on $\cC_k$
\be
f(y_k)\le f(y_{k^*})+\mu_{k^*}R_{k^*}^\rho.
\ee{afterk-s}
\EL
\proof
Note that for $k = 1$ \rf{rbeforek-s} is valid.
 Assume
it is valid for some $k \geq 1$. Note that
by \rf{NB1-mdp} of Corollary \rc{Stage1-mdp} there exists a random set, let us call it $\cB_k$, such that
$\Prob[\cB_k]\ge 1-\bar\al$ and on $\cB_k$,
\be
\delta_N(y_{k-1},R_{k-1})&\le& 2R_{k-1} \left(\sqrt{(L^2+\sigma^2) A(d) \over 2\mu(d) (N_0+1)}+
\sigma \sqrt{3\ln \bar \al^{-1}\over N_0+1}\right)\nn
&=&R_{k-1}a(N_0,\bar \alpha)\stackrel{\rf{mui}}{=}{1\over 2}\mu_k2^{-(\rho-1)k}R_0^{\rho-1}R_{k-1}
=\mu_kR^{\rho}_k.
\ee{ex36}
On the other hand, by our inductive hypothesis, $\| y_{k-1} - x^* \| \leq R_{k-1}$ on $\cA_{k-1}$.
Let
$\cA_k=\cA_{k-1}\cap \cB_k$. Note that
\[
\Prob[\cA_k]\ge \Prob[\cA_{k-1}]+\Prob[\cB_k]-1\ge 1-k\bar\al,
\]
and
we have on $A_k$:
\bse
f(y_k)-f^*&\le &\delta_N(y_{k-1},R)\le \mu_kR^{\rho}_k,\\
\| y_{k} - x^* \|^\rho&\le& {\delta_N(y_{k-1},R_{k-1})\over \mu(f)}\le R^{\rho}_k,
\ese
what is \rf{dlbeforek-s} and \rf{rbeforek-s} for  $k+1$.
\[
\]
To show \rf{afterk-s} notice that,
we have for $k>k^*$ (cf. \rf{ex36})
\[
f(y_{k}) - f(y_{k-1})\le \dl_{N_0}(y_{k-1},R_{k-1})\le \mu_kR_{k}^\rho
\]
on some $\cB_k\subset\Omega$ such that $\Prob[\cB_k]\ge 1-\bar\al$.
Then we have on $\cC_k=\cap_{j=k^*+1}^k \cB_j$:
\[
\ba{rcl}
f(y_{k}) - f(y_{k^*})& =&\sum\limits_{j=k^*+1}^k
f(y_{j}) - f(y_{j-1}) \le \sum\limits_{j=k^*+1}^{k}
2^{k^*-j}\mu_{k^*} R^\rho_{k^*} \le \mu_{k^*}R^\rho_{k^*}.
\ea
\]
Note that $\Prob[\cC_k]\ge 1-(k-k^*)\bar\al$. This proves the lemma.
\qed
\par
Now we can finish the proof of the theorem. Let
$\mu_0\le \mu(f)\le \mu_m$. At the end of the $k^*$-th
stage we have on the set $\cA_{k^*}$ of probability at least $1-k^*\bar\al$:
\bse
f({y}_{k^*}) - f^* & \leq &
\dl_{N_0}(y_{k^*-1},R_{k^*-1})\; {\leq} \;
 \mu_{k^*} R_{k^*}^\rho.
\ese
Then on the set $\cA_{k^*}\cap \cC_m$ such that $\Prob[\cA_{k^*}\cap \cC_m]\ge 1-m\bar\al$ (cf \rf{afterk-s}) we have
\bse
f({y}_m) - f^*&\le&  2\mu_{k^*} R_{k^*}^\rho\le
 4{\mu_{k^*}^{1\over \rho-1}\over \mu(f)^{1\over \rho-1}}\mu_{k^*}R^{\rho}_{k^*}
 \\
 &=&4{\mu_0^{\rho\over \rho-1}R_0^\rho\over \mu(f)^{1\over \rho-1}}\;\stackrel{\rf{mui}}{\le}
 4\,
\left({2 a(N_0,\bar\al)\over \mu(f)^{1\over \rho}}\right)^{\rho\over \rho-1}.
\ese
It suffices to recall now that by the definition of $m$, $m\le {1\over 2}\log_2N$,
thus $m\bar\al\le \al$.
\par
If $\mu(f)< \mu_0$, we have on $\cA_1=\cB_1$ (cf. \rf{ex36}):
\bse
f(y_1)-f^*&\le& R_0a(N_0,\bar\al)={R_0\over a(N_0,\bar\al)^{1\over \rho-1}}\,
a(N_0,\bar\al)^{\rho\over \rho-1}\\
&=&
2^{1\over \rho-1}{a(N_0,\bar\al)^{\rho\over \rho-1}\over \mu_0^{1\over \rho-1}}\le
2^{1\over \rho-1}\left({a(N_0,\bar\al)\over \mu(f)^{1\over \rho}}\right)^{\rho\over \rho-1}.
\ese
Finally, we conclude using \rf{afterk-s}:  on $\cA_1\cap \cC_m$ we have
\[
f(y_m)-f^*\le 2R_0a(N_0,\bar\al)\le \left({2a(N_0,\bar\al)\over \mu(f)^{1\over \rho}}\right)^{\rho\over \rho-1}.
\]
\qed
\end{document}